\documentclass[11pt,a4paper]{article}
\usepackage{multirow}
 \usepackage{epsfig}
 \usepackage{epstopdf}
\usepackage[T1]{fontenc}
\usepackage{geometry}
\usepackage{amsbsy,amsmath,latexsym,amsfonts, epsfig, color, authblk, amssymb, graphics, bm}
\usepackage{epsf,slidesec,epic,eepic}
\usepackage{fancybox}
\usepackage{fancyhdr}
\usepackage{setspace}
\usepackage{nccmath}
\usepackage[normalem]{ulem}
\usepackage[colorlinks, citecolor=blue]{hyperref}

\newtheorem{theorem}{Theorem}[section]
\newtheorem{lemma}{Lemma}[section]

\newtheorem{definition}{Definition}[section]
\newtheorem{example}{Example}[section]
\newtheorem{proposition}{Proposition}[section]
\newtheorem{corollary}{Corollary}[section]

\newtheorem{remark}{Remark}[section]

\newenvironment{proof}{{\noindent \bf Proof:}}{\hfill$\Box$\medskip}

\definecolor{lred}{rgb}{1,0.8,0.8}
\definecolor{lblue}{rgb}{0.8,0.8,1}
\definecolor{dred}{rgb}{0.6,0,0}
\definecolor{dblue}{rgb}{0,0,0.5}
\definecolor{dgreen}{rgb}{0,0.5,0.5}

\title{Explicit KL Exponents for PLQ Functions on the Sphere}
\author{Yanping Chen\footnote{School of Mathematics, South China University of Technology, Guangzhou, China.},\ \ Shaohua Pan\footnote{School of Mathematics, South China University of Technology, Guangzhou (shhpan@scut.edu.cn).}\ \ {\rm and}\ \ Anna Qi\footnote{School of Mathematics and Systems Science, Guangdong Polytechnic Normal University, Guangzhou (qiann@gpnu.edu.cn).}}

\begin{document}

\maketitle

\begin{abstract}
Let $\phi$ be a proper piecewise linear-quadratic (PLQ) function
whose domain meets the unit sphere $\mathcal{S}$. We prove that
$f:=\phi+\delta_{\mathcal{S}}$, where $\delta_{\mathcal{S}}$ is
the indicator function of $\mathcal{S}$, satisfies the Kurdyka--{\L}ojasiewicz (KL) property with exponent $3/4$ throughout its domain. This exponent is sharp even for the convex PLQ class. We characterize the smallest KL exponents of quadratics restricted to affine sections of $\mathcal{S}$, obtaining the values $0$, $1/2$, $2/3$, and $3/4$. A tangential projection identity shows that, at each point of ${\rm dom}\,f$, the maximum of the exponents from the relevant affine sections is admissible for $f$. When $\phi$ is piecewise affine (PA), each of these smallest exponents is $0$ or $1/2$. Thus $f$ satisfies the KL property with exponent $1/2$ throughout ${\rm dom}\,f$, and this exponent is sharp over the class of sphere-constrained PA functions. When ${\rm dom}\,\phi$ is a finite union of closed polyhedral cones, $f$ admits KL exponent $1/2$ for homogeneous quadratic pieces and exponent $0$ for nonnegative $\phi$ that is linear on each cone. For any convex polytope $C$ containing the origin, we use the KL exponent $0$ to show that the full-step manifold proximal point iterates for minimizing the support function of $C$ over $\mathcal{S}$ stabilize at a critical point in finitely many steps. 
\end{abstract}


\section{Introduction}\label{sec1}

Let $\mathbb{E}$ be a finite-dimensional Euclidean space with
inner product $\langle\cdot,\cdot\rangle$ and induced norm
$\|\cdot\|$, and let $\mathcal{S}:=\{z\in\mathbb{E}\mid\|z\|=1\}$
be its unit sphere. We study the Kurdyka--{\L}ojasiewicz (KL) property with explicit exponents for 
\begin{equation}\label{eq:sphere-function}
f:=\phi+\delta_{\mathcal{S}},
\end{equation}
where $\phi:\mathbb{E}\to\mathbb{R}\cup\{+\infty\}$ is a proper piecewise linear-quadratic (PLQ) function whose effective domain meets $\mathcal{S}$, and $\delta_{\mathcal{S}}$ is the indicator function of the unit sphere $\mathcal{S}$. 

Several sphere-constrained problems can be formulated by
minimizing $f$ over $\mathbb{E}$. Nonnegative principal component
analysis (PCA) maximizes a quadratic form over the
nonnegative unit sphere~\cite{Montanari2016}. Testing copositivity
amounts to checking whether a quadratic form has a nonnegative
minimum over the same set~\cite{Ferreira2024}. After negating the PCA objective, we can write both quadratic optimization problems in the form
$\min_{z\in\mathbb{E}}f(z)$ with $f$ defined in \eqref{eq:sphere-function}. In each case, $\phi$ is the sum of a quadratic form and the indicator function of the nonnegative orthant. For a given matrix $A\in\mathbb{R}^{m\times n}$, minimizing $\|Az\|_1$ over unit
vectors $z\in\mathbb{R}^n$ arises in dual principal component
pursuit~\cite{TsakirisVidal2018} and in sparse vector recovery
after parametrizing the subspace by an orthonormal
basis~\cite{QuSunWright2016}. For this minimization problem, Chen et al.~\cite{ChenSo2021} developed a manifold proximal point algorithm (ManPPA) and established local quadratic convergence in distance to a set of weak sharp minima for sufficiently small proximal parameters. Maximizing $\|Az\|_1$ over unit vectors gives the single-component $\ell_1$-norm PCA
problem~\cite{Kwak2008}. Both $\|Az\|_1$ and $-\|Az\|_1$ are
piecewise-affine (PA) since $\|Az\|_1=\max_{s\in\{-1,1\}^m}\langle A^{\top}s,z\rangle$. Other examples include sphere-constrained $\ell_1$-regularized quadratic minimization~\cite{Shen2022} and the Huber formulation of an oblique Procrustes 
subproblem~\cite[Eq.~(5)]{Trendafilov2004}.

The function $f$ is proper, lower semicontinuous (lsc), and semialgebraic, so it satisfies the KL property by \cite{Bolte2007-II}. This property ensures
that any sequence satisfying the sufficient decrease, relative
error, and subsequential continuity conditions converges to
a critical point of $f$ (see \cite[Theorem~2.9]{Attouch2013}). An explicit KL exponent at the critical point further provides rate estimates for both the objective values and the iterates~\cite[Theorem~3.4]{Frankel2015}. To obtain such exponents, Li and Pong \cite[Section~3]{LiPong2018} developed calculus rules for deriving KL exponents from known ones under operations such as finite pointwise minima and certain smooth compositions, and Yu et al.\ \cite[Theorem~3.1]{YuLiPong2022} extended the finite-minimum
rule to inf-projections under suitable assumptions. A proper PLQ function is continuous on its domain and can be written as $\min_i(q_i+\delta_{P_i})$, where $q_i$ is the quadratic piece on the closed polyhedron $P_i$. Applying \cite[Corollary~5.2]{LiPong2018} to this representation gives the KL exponent $1/2$ for $\phi$ in the ambient space. For the sphere-constrained function $f$, the corresponding representation is
$\min_i\bigl(q_i+\delta_{P_i\cap\mathcal{S}}\bigr)$. The intersections $P_i\cap\mathcal{S}$ can be nonpolyhedral,
so the polyhedrality assumption of \cite[Corollary~5.2]{LiPong2018} need not hold for such representation of $f$.  

Explicit KL exponents are available for several classes of
functions restricted to the sphere. Liu et al.~\cite{LiuSo2019} established exponent $1/2$ for a class of homogeneous quadratic functions restricted to the Stiefel manifold, which reduces to the unit sphere in the single-column case. Wu et al.~\cite{WuPanBi2021} obtained exponent $1/2$ for zero-norm regularized homogeneous quadratic functions restricted to the unit sphere. For $\|Az\|_1+\delta_{\mathcal{S}}(z)$, Yu et al.\
\cite[Theorems~4.3 and~4.7]{YuZengPong2025} established
exponent $1/2$ when the optimal value is positive and exponent
$0$ at points in a set of weak sharp minima. For the trust-region subproblem over the closed unit ball, Jiang and Li~\cite[Theorem~4.8]{JiangLi2022} establish local KL
inequalities at global minimizers, with exponent $3/4$ in a
specific degenerate case and $1/2$ otherwise. Jiang and Li 
\cite[p.~1617]{JiangDuanLi2019} provided an example in which the
pointwise maximum of two convex quadratics has smallest KL
exponent $3/4$ at a minimizer where both quadratics are active.
These results motivate determining a sharp KL exponent valid
throughout ${\rm dom}\,f$ for general PLQ $\phi$ and deriving
pointwise refinements from the smallest KL exponents of the
quadratic pieces of $\phi$ restricted to the corresponding
affine sections of $\mathcal{S}$. 
\subsection{Main contributions}\label{sec1.1}

We establish in Theorem~\ref{theorem:PLQ-34KL}(ii) that $f$ satisfies the KL property with exponent $3/4$ at every point of ${\rm dom}\,f$, independently of the ambient dimension. Example~\ref{exa:PLQ-sharp} shows that this exponent is sharp over the PLQ class by exhibiting a nonsmooth interface point at which the smallest KL exponent of $f$ is $3/4$. At that
point, each adjacent quadratic piece has a nonzero tangential
gradient on $\mathcal{S}$ and therefore has smallest KL
exponent $0$ when restricted to the full sphere.

Theorem~\ref{theorem:PLQ-34KL}(i) provides an admissible KL
exponent $\theta_{\overline{z}}$, possibly smaller than $3/4$,
for $f$ at each $\overline{z}\in{\rm dom}\,f$. To obtain $\theta_{\overline{z}}$, we consider the cells
$C\subseteq{\rm dom}\,\phi$ whose closures contain $\overline{z}$.
For each such cell, Proposition~\ref{prop:quadratic-sphere-KL}
characterizes the smallest KL exponent at $\overline{z}$ of the
corresponding quadratic piece restricted to
${\rm aff}(C)\cap\mathcal{S}$. The possible values are $0$,
$1/2$, $2/3$, and $3/4$. Using the tangential projection identity in
Proposition~\ref{prop:tangent-proj-PLQ}, we combine the
gradient estimates on these sections to obtain a KL inequality
for $f$ at $\overline{z}$ with exponent $\theta_{\overline{z}}$
equal to the maximum of these smallest KL exponents. 

For proper PA $\phi$, Theorem~\ref{theorem:KL-PA} establishes KL exponent $1/2$ throughout ${\rm dom}\,f$, and Remark~\ref{remark-PA1/2} shows that this exponent is sharp over the PA class. Proposition~\ref{prop:ph12} establishes the same exponent when ${\rm dom}\,\phi$ is a finite union of closed convex polyhedral cones and $\phi$ is homogeneous quadratic on each
cone. In this case, for each reference point
$\overline{z}\in{\rm dom}\,f$, the spectral estimate applies
at every $z\in{\rm dom}\,f$ with a sufficiently small positive
function-value gap $f(z)-f(\overline{z})$. 
Example~\ref{exa:PLQ-sharp} also shows that the sum of a PA
function and a homogeneous quadratic function, both restricted
to $\mathcal{S}$, can have smallest KL exponent $3/4$ at a
point where each summand admits exponent $1/2$.
If $\phi$ is nonnegative and linear on each of finitely many
closed polyhedral cones whose union is ${\rm dom}\,\phi$,
Proposition~\ref{prop:KL-PL} gives exponent $0$ throughout
${\rm dom}\,f$. 

For any convex polytope $C$ containing the origin, Corollary~\ref{cor:manppa-finite} uses this KL exponent $0$ to show that full-step ManPPA for minimizing the support function $\sigma_C$ over $\mathcal{S}$ stabilizes at a critical point in finitely many steps. The guarantee holds for every initial point on $\mathcal{S}$ and every fixed positive proximal
parameter. For $C=A^\top[-1,1]^m$, so that $\sigma_C(z)=\|Az\|_1$,
the full-step iteration coincides with \cite[Algorithm~1]{ChenSo2021} when the tangent subproblems are solved exactly, since the line search accepts the unit step at every iteration. Our finite convergence guarantee strengthens the conclusion of \cite[Theorem~1]{ChenSo2021} that every
cluster point is stationary. For this fixed-parameter 
$\ell_1$ iteration, it also removes the assumption in the
finite convergence theorem of Yu et al.~\cite{YuZengPong2025}
that the limit lies in a set of weak sharp minima.

\subsection{Notation}\label{sec1.2}

Throughout this paper, set $\overline{\mathbb{R}}:=\mathbb{R}\cup\{+\infty\}$ and let $\mathcal{I}$ denote the identity map on $\mathbb{E}$. For each positive integer $k$, set $[k]:=\{1,\ldots,k\}$. For $x\in\mathbb{E}$ and $D\subseteq\mathbb{E}$, define
${\rm dist}(x,D):=\inf_{y\in D}\|x-y\|$, with the convention $\inf\varnothing:=+\infty$. For a set $C\subseteq\mathbb{E}$,  ${\rm cl}\,C$, ${\rm aff}(C)$, and ${\rm ri}(C)$ denote the closure, affine hull, and relative interior of $C$, respectively. The indicator function
$\delta_C$ is defined by $\delta_C(x)=0$ if $x\in C$ and
$\delta_C(x)=+\infty$ otherwise. If $C$ is closed and $x\in C$, then  $\mathcal{N}_C(x)$ and $T_x C$ denote the limiting normal cone and the
(Bouligand) tangent cone to $C$ at $x$, respectively. For $x\in\mathbb{E}$ and $\varepsilon>0$, let $\mathbb{B}(x,\varepsilon)$ denote the closed ball centered at $x$ with radius $\varepsilon$, and let $\mathbb{B}^{\circ}(x,\varepsilon)$ denote its interior. For a linear subspace $\mathbb{L}\subseteq\mathbb{E}$, let $\mathcal{P}_{\mathbb{L}}$ be the orthogonal projection onto $\mathbb{L}$ and 
$\mathcal{I}_{\mathbb{L}}$ the identity map on $\mathbb{L}$.
Set $\mathbb{B}_{\mathbb{L}}(0,\varepsilon)
:=\mathbb{L}\cap\mathbb{B}(0,\varepsilon)$ and $\mathbb{B}_{\mathbb{L}}^{\circ}(0,\varepsilon)
:=\mathbb{L}\cap\mathbb{B}^{\circ}(0,\varepsilon)$ for $\varepsilon>0$. For a proper lsc function $h:\mathbb{E}\to\overline{\mathbb{R}}$ and a point $x\in{\rm dom}\,h$, $\widehat{\partial}h(x)$ and $\partial h(x)$ denote the regular and limiting subdifferential of $h$ at $x$, respectively. For real-valued functions $a$ and $b$ defined near $t_0$, we write $a(t)\sim b(t)$
as $t\to t_0$ if $b(t)\ne0$ for all $t\ne t_0$ sufficiently
close to $t_0$ and $\lim_{t\to t_0}{a(t)}/{b(t)}=1$. 
\section{Preliminaries}\label{sec2}

We give the formal definitions of PLQ, PA, and PL functions.
\begin{definition}\label{PLQ-def}(see \cite[Definition 10.20]{RW98})
 A function $h:\mathbb{E}\to\overline{\mathbb{R}}$ is called piecewise linear-quadratic if ${\rm dom}\,h$ can be represented as the union of finitely many polyhedral sets $P_1,\ldots,P_m$, relative to each of which $h$ is given by $h(x)=\frac12\langle x,\mathcal{Q}_ix\rangle+\langle a_i,x\rangle+b_i$ for some self-adjoint linear operator $\mathcal{Q}_i:\mathbb{E}\to\mathbb{E}$, vector $a_i\in\mathbb{E}$, and scalar $b_i\in\mathbb{R}$. When $\mathcal{Q}_i\equiv 0$ for all $i\in[m]$, it is called piecewise-affine; when, in addition, $b_i\equiv 0$ for all $i\in[m]$, it is called piecewise-linear.
\end{definition}

By \cite[Proposition 10.21]{RW98}, a PLQ function $h:\mathbb{E}\to\overline{\mathbb{R}}$ has a closed domain ${\rm dom}\,h$ and is continuous relative to ${\rm dom}\,h$, hence lower semicontinuous (lsc) on $\mathbb{E}$. 

We next recall the KL property with exponent $\theta\in[0,1)$ for a nonsmooth function.
\begin{definition}\label{def-KL}(see \cite{Bolte2007,LiPong2018})
A proper lsc function $h\!:\mathbb{E}\to\overline{\mathbb{R}}$ is said to have the KL property with exponent $\theta\in[0,1)$ at $\overline{x}\in{\rm dom}\,\partial h$ if there exist $c>0,\,\varepsilon>0$ and $\eta\in(0,\infty]$ such that, for all $x\in\mathbb{B}(\overline{x},\varepsilon)\cap[h(\overline{x})<h<h(\overline{x})+\eta]$, 
\[
 {\rm dist}(0,\partial h(x))\ge c(h(x)-h(\overline{x}))^{\theta}.
\]
If $h$ has the KL property with the same exponent $\theta$ at every point of ${\rm dom}\,\partial h$, it is called a KL function of exponent $\theta$.
\end{definition}

For a proper lsc semialgebraic function $h$, the nonsmooth KL
inequality in \cite[Theorem~14]{Bolte2007-II}
and the growth dichotomy for one-variable semialgebraic
functions yield a KL exponent at each point of
${\rm dom}\,\partial h$. Obtaining an explicit admissible
exponent for a given class generally requires quantitative
estimates tailored to its structure.

In the remainder of this section, we establish several auxiliary results for use in the subsequent analysis. We begin by examining the incidence relations among the relatively open polyhedral cells induced by a finite hyperplane arrangement. The following lemma shows that if one such cell intersects the closure of another, then it is contained in that closure. It also establishes the corresponding inclusion relations between their affine hulls and between the linear subspaces parallel to those hulls. These geometric properties will be instrumental in establishing the consistency of the quadratic representations of PLQ functions over the induced cells. 
\begin{lemma}\label{lem:hplane-partition}
 Let $H_1,\ldots,H_q$ be finitely many hyperplanes in $\mathbb{E}$, and let $C_1,\ldots,C_N$ denote all the nonempty relatively open polyhedral cells obtained by partitioning $\mathbb{E}$ according to the signs of the affine functions defining these hyperplanes, including all lower-dimensional cells. For each $j\in[N]$, define $\mathbb{A}_j:={\rm aff}(C_j)$ and $\mathbb{L}_j:=\mathbb{A}_j-\mathbb{A}_j$. Then, for any $j,\ell\in[N]$, 
 \begin{itemize}
 \item[(i)] if $C_j\cap{\rm cl}\,C_\ell\neq\varnothing$, then $C_j\subseteq{\rm cl}\,C_\ell$; 
 
 \item[(ii)] if $C_j\cap{\rm cl}\,C_\ell\neq\varnothing$, then $\mathbb{A}_j		\subseteq\mathbb{A}_\ell$ and $\mathbb{L}_j\subseteq\mathbb{L}_{\ell}$. 
 \end{itemize}
 \end{lemma}
\begin{proof}
If $q=0$, the empty arrangement induces the single cell
$C_1=\mathbb{E}$, with $\mathbb{A}_1=\mathbb{L}_1=\mathbb{E}$,
so (i) and (ii) hold immediately. Henceforth we assume $q\ge1$. 
For each $s\in[q]$, let $h_s:\mathbb{E}\to\mathbb{R}$ be a nonconstant affine function defining the hyperplane $H_s$, namely, $H_s=\{z\in\mathbb{E}\mid h_s(z)=0\}$. For any $x\in\mathbb{E}$, define 
\[
 \sigma(x):=\big({\rm sgn}(h_1(x)),\ldots,{\rm sgn}(h_q(x))\big)^{\top}\ \ {\rm and}\ \ \Sigma:=\big\{\sigma(x)\mid x\in\mathbb{E}\big\},
\]
where ${\rm sgn}(t)=1$ if $t>0$, ${\rm sgn}(t)=-1$ if $t<0$, and ${\rm sgn}(t)=0$ if $t=0$. For each $\xi\in\Sigma$, let
$C_{\xi}:=\{x\in\mathbb{E}\mid \sigma(x)=\xi\}$. By the definition of $\Sigma$, each $C_\xi$ is nonempty, and the sets $C_\xi$ for $\xi\in\Sigma$ are pairwise disjoint and have union $\mathbb{E}$.
They are precisely the nonempty relatively open polyhedral cells induced by $H_1,\ldots,H_q$. Since $\Sigma\subseteq\{-1,0,1\}^q$ is finite and its elements are in one-to-one
correspondence with these cells, we have $N=|\Sigma|$. Thus, we may enumerate the elements of $\Sigma$ as
$\sigma^1,\ldots,\sigma^N$, where
$\sigma^j=(\sigma_1^j,\ldots,\sigma_q^j)^\top$.
Accordingly, we may label the cells so that,
for each $j\in[N]$,
\[
 C_j=\Big\{x\in\mathbb{E}\mid\sigma(x)=\sigma^j\Big\}=\Bigg\{
 x\in\mathbb{E}\,\bigg|
 \begin{array}{ll}
 h_s(x)>0 & \text{if }\sigma_s^j=1,\\
 h_s(x)=0 & \text{if }\sigma_s^j=0,\\
 h_s(x)<0 & \text{if }\sigma_s^j=-1
 \end{array}\ {\rm for}\ s=1,\ldots,q\Bigg\}.
\]
For each $j\in[N]$, define $\mathcal H_j^0
:=\bigcap_{\substack{s\in[q],\sigma_s^j=0}}H_s$, $\mathcal O_j^+
:=\bigcap_{\substack{s\in[q],\sigma_s^j=1}}
\{x\in\mathbb{E}\mid h_s(x)>0\}$, and $\mathcal O_j^-
:=\bigcap_{\substack{s\in[q],\sigma_s^j=-1}}
\{x\in\mathbb{E}\mid h_s(x)<0\}$, where, if the corresponding index set is empty, the associated set is understood to be $\mathbb E$. Then  $C_j=\mathcal{H}_j^0\cap\mathcal{O}_j^+\cap\mathcal{O}_j^-$. Since $C_j\neq\varnothing$, $\mathcal{H}_j^0$ is a nonempty affine subspace of
$\mathbb{E}$. Moreover, since $\mathcal{O}_j^+$ and $\mathcal{O}_j^-$ are open subsets of $\mathbb{E}$, $C_j$ is a relatively open subset of $\mathcal{H}_j^0$. Consequently, $\mathbb A_j={\rm aff}(C_j)=\mathcal H_j^0$, and $C_j$ is relatively open in $\mathbb A_j$. By construction, $C_1,\ldots,C_N$ are nonempty, pairwise disjoint, relatively open polyhedral cells satisfying $\mathbb{E}=\bigcup_{j=1}^N C_j$. Furthermore, for each $j\in[N]$,  \begin{equation}\label{equa-closure}
  {\rm cl}\,C_j=\Bigg\{z\in\mathbb{E}\,\bigg|
  \begin{array}{ll}
  h_s(z)\geq0 &\text{if }\sigma_s^j=1,\\
  h_s(z)=0 &\text{if }\sigma_s^j=0,\\
  h_s(z)\leq0 &\text{if }\sigma_s^j=-1
 \end{array}\ {\rm for}\ s=1,\ldots,q\Bigg\}.
 \end{equation}
 Indeed, let $D_j$ denote the set on the right hand side of \eqref{equa-closure}. The inclusion ${\rm cl}\,C_j\subseteq D_j$ follows by the continuity of the affine functions $h_1,\ldots,h_q$. Conversely, fix any $\widehat{z}\in C_j$ and take an arbitrary $z\in D_j$. For any $t\in(0,1]$, let $z_t:=(1-t)z+t\widehat z$. Then we have $h_s(z_t)=(1-t)h_s(z)+t h_s(\widehat z)$ for each $s\in[q]$. By the definitions of $D_j$ and $C_j$, 
 \[
	h_s(z_t)>0\ \ \text{if }\sigma_s^j=1,\ h_s(z_t)=0\ \ \text{if }\sigma_s^j=0,\ \ {\rm and}\ \ h_s(z_t)<0\ \ \text{if }\sigma_s^j=-1.
 \]
 Consequently, $z_t\in C_j$ for every $t\in(0,1]$. Since $z_t\to z$ as $t\downarrow0$, we obtain $z\in{\rm cl}\,C_j$. The arbitrariness of $z\in D_j$ implies that $D_j\subseteq{\rm cl}\,C_j$, and hence \eqref{equa-closure} holds.

 \noindent
 {\bf(i)} Suppose that $C_j\cap{\rm cl}\,C_{\ell}\neq\varnothing$, and take any $z\in C_j\cap{\rm cl}\,C_{\ell}$. For any $u\in C_j$, since $u$ and $z$ belong to the same cell, we have $\sigma(u)=\sigma(z)$. On the other hand, since $z\in{\rm cl}\,C_{\ell}$, 	it satisfies all the weak inequalities and equalities in the characterization \eqref{equa-closure} of ${\rm cl}\,C_{\ell}$. Since $\sigma(u)=\sigma(z)$, the signs of $h_s(u)$ and $h_s(z)$ coincide for every $s\in[q]$. Consequently, $u$ also satisfies all these weak inequalities and equalities, which implies that $u\in{\rm cl}\,C_\ell$. Since $u\in C_j$ is arbitrary, we obtain $C_j\subseteq{\rm cl}\,C_\ell$. This establishes (i). 
 
 \noindent
 {\bf(ii)} Under the same assumption, item (i) gives
$C_j\subseteq{\rm cl}\,C_\ell$. Since $\mathbb A_\ell$ is
a closed affine subspace containing $C_\ell$, we have
${\rm cl}\,C_\ell\subseteq\mathbb{A}_\ell$. It follows that
$\mathbb{A}_j={\rm aff}(C_j)\subseteq\mathbb{A}_\ell$. Hence
$\mathbb{L}_j=\mathbb{A}_j-\mathbb{A}_j\subseteq \mathbb{A}_{\ell}-\mathbb{A}_{\ell}=\mathbb{L}_{\ell}$. Thus, we complete the proof.
\end{proof}

Combining Lemma~\ref{lem:hplane-partition} with the piecewise quadratic representation of a PLQ function, we obtain the following compatibility properties of its quadratic pieces.
\begin{corollary}\label{cor:PLQ-consistency}
Let $\phi:\mathbb{E}\to\overline{\mathbb{R}}$ be a proper PLQ function with ${\rm dom}\,\phi=\bigcup_{i=1}^m P_i$, where $P_1,\ldots,P_m$ are closed polyhedral sets, and $\phi(z)=\frac12\langle z,\mathcal{Q}_i z\rangle+\langle a_i,z\rangle+b_i$ for $z\in P_i$. Here, $\mathcal{Q}_i:\mathbb{E}\to\mathbb{E}$ is a self-adjoint linear operator, $a_i\in\mathbb{E}$ and $b_i\in\mathbb{R}$. Collect all the hyperplanes appearing in the halfspace and hyperplane representations of $P_1,\ldots,P_m$, and denote them by $H_1,\ldots,H_q$. Let $C_1,\ldots,C_N$ and, for each $j\in[N]$, $\mathbb A_j$ and $\mathbb L_j$ be as defined in Lemma~\ref{lem:hplane-partition}. Then, the following assertions hold. 
 \begin{itemize}
 \item[(i)] For any $i\in[m]$ and $j\in[N]$, either $C_j\subseteq P_i$ or $C_j\cap P_i=\varnothing$.
  	
 \item[(ii)] For each $j$ such that $C_j\subseteq{\rm dom}\,\phi$, there exists an index $i_j\in[m]$ such that $C_j\subseteq P_{i_j}$. For each such $j$, fix one such index $i_j$ and set
 $(\widetilde{\mathcal{Q}}_j,\widetilde{a}_j,\widetilde{b}_j):=(\mathcal{Q}_{i_j},a_{i_j},b_{i_j})$. 
 
 \item[(iii)] For any $j,\ell$ such that $C_j\cup C_{\ell}\subseteq{\rm dom}\,\phi$ and $C_j\cap{\rm cl}\,C_\ell\neq\varnothing$, it holds
 \[
 \langle (\widetilde{\mathcal{Q}}_j-\widetilde{\mathcal{Q}}_\ell)x+\widetilde{a}_j-\widetilde{a}_\ell,d\rangle=0\quad\forall x\in C_j,\,\forall d\in\mathbb{L}_j.
 \]
 \end{itemize}  
\end{corollary}
\begin{proof}
{\bf(i)} Fix any $i\in[m]$ and $j\in[N]$. Suppose that $C_j\cap P_i\neq\varnothing$. We prove that $C_j\subseteq P_i$. Take any $\widehat{x}\in C_j\cap P_i$. Clearly, the point $\widehat{x}$ satisfies all the defining inequalities and equalities of $P_i$. Since each defining affine function has the same sign at every point of $C_j$ as at $\widehat{x}$, these constraints hold throughout $C_j$.
Hence, $C_j\subseteq P_i$, and item (i) follows. 
 
\noindent
{\bf(ii)} Fix any $j\in[N]$ such that $C_j\subseteq{\rm dom}\,\phi$. Since $C_j\ne\varnothing$ by construction, we may choose $\widehat{x}^j\in C_j$. By
${\rm dom}\,\phi=\bigcup_{i=1}^m P_i$ and $C_j\subseteq{\rm dom}\,\phi$, there exists an index $i_j\in[m]$ such that $\widehat{x}^j\in P_{i_j}$. Consequently, $C_j\cap P_{i_j}\neq\varnothing$. By item (i), we obtain $C_j\subseteq P_{i_j}$. 
 
 \noindent
 {\bf(iii)}
Fix any $j,\ell$ satisfying the stated conditions. By Lemma~\ref{lem:hplane-partition}(i), $C_j\subseteq{\rm cl}\,C_{\ell}$. Since $C_{\ell}\subseteq P_{i_{\ell}}$ and $P_{i_{\ell}}$ is closed, we have ${\rm cl}\,C_\ell\subseteq P_{i_\ell}$. Together with
$C_j\subseteq P_{i_j}$ by item (ii), this gives $C_j\subseteq P_{i_j}\cap P_{i_\ell}$. Thus, for every $z\in C_j$, 
 \[
 \frac12\langle z,\widetilde{\mathcal{Q}}_jz\rangle+\langle \widetilde{a}_j, z\rangle+\widetilde{b}_j=\phi(z)=\frac12\langle z,\widetilde{\mathcal{Q}}_\ell z\rangle+\langle \widetilde{a}_\ell,z\rangle+\widetilde{b}_{\ell}.
 \]
 Fix any $x\in C_j$ and $d\in\mathbb{L}_j$. Since
$\mathbb{A}_j=x+\mathbb{L}_j$ and $C_j$ is open relative to
$\mathbb{A}_j$, there exists $\varepsilon>0$ such that
$x+td\in C_j$ for all $t\in(-\varepsilon,\varepsilon)$.
Substituting $z=x+td$ and differentiating the two quadratic
expressions at $t=0$ gives the desired identity.
\end{proof}
\section{Main results}\label{sec3}

We establish that the function $f$ defined in \eqref{eq:sphere-function}
satisfies the KL property throughout ${\rm dom}\,f$ with exponent
$3/4$ when $\phi$ is PLQ and with exponent $1/2$ when $\phi$ is PA.
Both exponents are sharp over their respective classes. Since the
KL property is defined at points in ${\rm dom}\,\partial f$, we
first show that ${\rm dom}\,\partial f={\rm dom}\,f$.
\begin{lemma}\label{lem:domain-Lsubdiff}
Let $\phi$ and its PLQ representation be as in Corollary~\ref{cor:PLQ-consistency}, and suppose that ${\rm dom}\,\phi\cap\mathcal{S}\ne\varnothing$. Then ${\rm dom}\,\partial f={\rm dom}\,f$.
\end{lemma}
\begin{proof}
The inclusion ${\rm dom}\,\partial f\subseteq{\rm dom}\,f$ is trivial. It remains to establish the converse inclusion. Fix any $\overline{x}\in{\rm dom}\,f$. Since ${\rm dom}\,\phi=\bigcup_{i=1}^m P_i$ and $\mathcal S$ are closed, so is
${\rm dom}\,f={\rm dom}\,\phi\cap\mathcal S$. Define the index set  
$I(\overline{x}):=\{i\in[m]\mid \overline{x}\in P_i\}$. Clearly, $I(\overline{x})\neq\varnothing$. If $I(\overline{x})=[m]$, choose any $\varepsilon>0$. Otherwise, the intersection
$\bigcap_{i\in[m]\setminus I(\overline{x})}(\mathbb{E}\setminus P_i)$ 
is an open set containing $\overline{x}$, so there exists
$\varepsilon>0$ such that $P_i\cap\{x\in\mathbb E\mid\|x-\overline{x}\|<\varepsilon\}=\varnothing$ for all $i\notin I(\overline{x})$. In either case, every $x\in{\rm dom}\,\phi\cap\mathbb{B}^{\circ}(\overline{x},\varepsilon)$ belongs to some $P_i$ with $i\in I(\overline{x})$. For each $i\in I(\overline{x})$, define the quadratic function $q_i:\mathbb{E}\to\mathbb{R}$ by
\[
 q_i(x):=\frac12\langle x,\mathcal{Q}_i x\rangle +\langle a_i,x\rangle+b_i. 
\]
Since each $q_i$ is continuously differentiable and $I(\overline{x})$ is finite, there exists a constant $L>0$ such that $\|\nabla q_i(x)\|\le L$ for all $i\in I(\overline{x})$ and $x\in\mathbb{B}(\overline{x},\varepsilon)$. For every $x\in{\rm dom}\,\phi\cap\mathbb{B}^{\circ}(\overline{x},\varepsilon)$, let $i\in I(\overline{x})$ be such that $x\in P_i$. Then, 
$\phi(x)=q_i(x)$ and $\phi(\overline{x})=q_i(\overline{x})$.
Applying the mean value theorem to $q_i$ along the segment
joining $\overline{x}$ and $x$ yields 
\begin{equation}\label{RLip-phi}
 |\phi(x)-\phi(\overline{x})|\le L\|x-\overline{x}\|\quad\text{for\ all}\ x\in {\rm dom}\,\phi\cap\mathbb{B}^{\circ}(\overline{x},\varepsilon).
\end{equation}
Since $\phi$ is continuous on ${\rm dom}\,\phi$ and ${\rm dom}\,f$ is nonempty and compact, for each $\lambda>0$, we can choose
\[
 x^\lambda\in\mathop{\arg\min}_{x\in{\rm dom}\,f}\Big\{ \phi(x)+\frac{\lambda}{2}\|x-\overline{x}\|^2\Big\}. 
\]
By the optimality of $x^{\lambda}$ and the feasibility
of $\overline{x}$, we have $\phi(x^\lambda)+\frac{\lambda}{2}\|x^\lambda-\overline{x}\|^2\le\phi(\overline{x})$. Since $\phi$ is bounded below on ${\rm dom}\,f$, this implies that
$x^\lambda\to\overline{x}$ as $\lambda\to\infty$. Therefore, for all sufficiently large $\lambda$, it holds $\|x^\lambda-\overline{x}\|<\varepsilon$. From \eqref{RLip-phi}, it then follows that $\frac{\lambda}{2}\|x^\lambda-\overline{x}\|^2\le \phi(\overline{x})-\phi(x^\lambda)\le L\|x^\lambda-\overline{x}\|$. Consequently, for all sufficiently large $\lambda$, 
\begin{equation}\label{temp-xlambda}
\lambda\|x^\lambda-\overline{x}\|\le 2L. 
\end{equation}
Note that, for each $\lambda>0$, $x^\lambda$ is also a global minimizer of
$f(\cdot)+\frac{\lambda}{2}\|\cdot-\overline{x}\|^2$ over $\mathbb{E}$. By the Fermat rule and \cite[Exercise 8.8(c)]{RW98}, we have $v^\lambda:=\lambda(\overline{x}-x^\lambda)\in\widehat{\partial} f(x^\lambda)$ for every $\lambda>0$. By \eqref{temp-xlambda}, $\|v^\lambda\|\le2L$ for all sufficiently large $\lambda$. Choose a sequence $\{\lambda_k\}\subset(0,\infty)$ such that $\lambda_k\to\infty$. Then $\{v^{\lambda_k}\}_{k\in\mathbb{N}}$ is bounded. Passing to a subsequence
and relabeling it, we may assume that $v^{\lambda_k}\to\overline{v}$
for some $\overline{v}\in\mathbb{E}$. Since $x^{\lambda_k}\to\overline{x}$, \eqref{RLip-phi} yields  $f(x^{\lambda_k})=\phi(x^{\lambda_k})\to\phi(\overline{x})=f(\overline{x})$. Together with $v^{\lambda_k}\in\widehat{\partial}f(x^{\lambda_k})$,
these limits give $\overline{v}\in\partial f(\overline{x})$ by the definition of the limiting subdifferential. Thus $\overline{x}\in{\rm dom}\,\partial f$. 
\end{proof}

We next show that the projection of each limiting subgradient
of $f$ onto the tangent space of the associated affine section
of $\mathcal S$ equals the gradient of the corresponding
quadratic piece of $\phi$ along that section. This 
allows us to use gradient estimates on these sections to
establish the KL property of $f$ with explicit exponents.
\begin{proposition}\label{prop:tangent-proj-PLQ}
Let $\phi$, its PLQ representation, and the cells
$C_1,\ldots,C_N$ be as in Corollary~\ref{cor:PLQ-consistency}.
Define $\mathcal{J}:=\{j\in[N]\mid C_j\subseteq{\rm dom}\,\phi\}$ and $\mathcal{M}_j:=\mathbb{A}_j\cap\mathcal{S}$ for each $j\in[N]$. For each $j\in\mathcal{J}$, let $i_j$ and
$(\widetilde{\mathcal{Q}}_j,\widetilde{a}_j,\widetilde{b}_j)$ be as in Corollary~\ref{cor:PLQ-consistency}(ii), and define $q_j(z):=\tfrac12\langle z,\widetilde{\mathcal{Q}}_j z\rangle +\langle \widetilde{a}_j,z\rangle+\widetilde{b}_j$ for $z\in\mathbb{E}$. Then, for every $j\in\mathcal{J}$ and $z\in C_j\cap\mathcal{S}$, 
\[ 
 \mathcal{P}_{T_z\mathcal{M}_j}v=\mathcal{P}_{T_z\mathcal{M}_j}\nabla q_j(z)\quad\ \forall v\in\partial f(z).
\]
\end{proposition}
\begin{proof}
 Fix any $j\in\mathcal{J}$ and $z\in C_j\cap\mathcal{S}$. We first establish the following equality
 \begin{equation}\label{Fsubdiff-PLQ}
\mathcal{P}_{T_z\mathcal{M}_j}v=\mathcal{P}_{T_z\mathcal{M}_j}(\widetilde{\mathcal{Q}}_jz+\widetilde{a}_j)\quad\ \forall v\in\widehat{\partial}f(z).
\end{equation}	 
 Let $p_j$ be the orthogonal projection of the origin onto $\mathbb{A}_j$. Then $\mathbb{A}_j=p_j+\mathbb{L}_j$  and $p_j\perp \mathbb{L}_j$. If $z=p_j$, then $\mathcal{M}_j=\{p_j\}=\{z\}$. 
 Indeed, for any $y\in\mathcal {M}_j$, by writing $y=p_j+\xi$ with $\xi\in\mathbb{L}_j$, we have
 $1=\|y\|^2=\|p_j\|^2+\|\xi\|^2=1+\|\xi\|^2$. Then $\xi=0$ and $y=p_j$, which implies the desired equalities. In this case, $T_z\mathcal M_j=\{0\}$ and \eqref{Fsubdiff-PLQ} holds trivially. Hence it suffices to assume that $z\neq p_j$. Let $r_j:=\|z-p_j\|=\sqrt{1-\|p_j\|^2}>0$. Since $\mathbb A_j=p_j+\mathbb L_j$ and $p_j\perp\mathbb L_j$, we have
 $\mathcal{M}_j=\{p_j+r_j u\mid u\in\mathbb{L}_j,\|u\|=1\}$.  Consequently, $T_z\mathcal{M}_j=\mathbb{L}_j\cap z^\perp$. If $T_z\mathcal{M}_j=\{0\}$, the desired equality \eqref{Fsubdiff-PLQ} holds. Otherwise, take any $0\ne d\in T_z\mathcal{M}_j$. Then, $d\perp z$ and $d\in\mathbb{L}_j$. Since $p_j\perp\mathbb{L}_j$, we have $d\perp p_j$. Set $u:=(z-p_j)/r_j$. Clearly, $u\in\mathbb{L}_j\cap\mathcal{S}$ and $d\perp u$. Consider the curve
 \[
	\gamma(t):=p_j+r_j\left(u\cos\frac{\|d\|t}{r_j}
	+\frac{d}{\|d\|}\sin\frac{\|d\|t}{r_j}\right)\quad{\rm for}\ t\in\mathbb{R}.
 \]
 Since $u,d\in\mathbb{L}_j$ and $p_j\perp\mathbb{L}_j$, we have $\gamma(t)\in\mathbb{A}_j$ and  $\|\gamma(t)\|=1$ for all $t\in\mathbb{R}$. Moreover, $\gamma(0)=p_j+r_ju=z$ and $\gamma'(0)=d$. Since $C_j$ is relatively open in $\mathbb{A}_j$ and $\gamma(0)=z\in C_j$, by the continuity of $\gamma$, there exists sufficiently small $\epsilon>0$ such that $\gamma(t)\in C_j\cap\mathcal{S}$ for all $t\in(-\epsilon,\epsilon)$. Since $C_j\subseteq P_{i_j}$, we have $f(\gamma(t))=\phi(\gamma(t))=q_j(\gamma(t))$ for all sufficiently small $|t|$, which together with $z\in C_j\cap\mathcal{S}$ and Corollary~\ref{cor:PLQ-consistency}(ii) yields 
 \[
	f(\gamma(t))-f(z)=\langle\widetilde{\mathcal{Q}}_jz+\widetilde{a}_j,\gamma(t)-z\rangle+o(\|\gamma(t)-z\|).
 \]
 On the other hand, since $\gamma'(0)=d\ne0$, we have $\gamma(t)-z=td+o(|t|)$
and $\|\gamma(t)-z\|/|t|\to\|d\|>0$ as $t\to0$.
Thus $\gamma(t)\ne z$ for all sufficiently small $t\ne0$. Consequently, for any $v\in\widehat{\partial}f(z)$, the definition of the Fr\'echet subdifferential gives 
\[
  \liminf_{t\to0}\frac{f(\gamma(t))-f(z)-\langle v,\gamma(t)-z\rangle}{\|\gamma(t)-z\|}\geq0.
\]
 Using $\gamma(t)-z=td+o(|t|)$, and letting $t\downarrow0$ and $t\uparrow0$, respectively, we obtain 
 \[
   \langle\widetilde{\mathcal{Q}}_jz+\widetilde{a}_j-v,d\rangle\ge 0\ \ {\rm and}\ \ \langle\widetilde{\mathcal{Q}}_jz+\widetilde{a}_j-v,d\rangle\le 0.
 \]
 Hence, $\langle v,d\rangle=\langle\widetilde{\mathcal{Q}}_jz+\widetilde{a}_j,d\rangle$, and \eqref{Fsubdiff-PLQ} follows by the arbitrariness of $d\in T_z\mathcal{M}_j$. 

We next establish the desired equality for the limiting subdifferential. Fix any $v\in\partial f(z)$. By the definition of the limiting subdifferential, there exist sequences $z^k\to z$ and $v^k\to v$ such that $f(z^k)\to f(z)$ and $v^k\in\widehat{\partial}f(z^k)$. For each $k\in\mathbb{N}$, since $v^k\in\widehat{\partial}f(z^k)$, we have $z^k\in{\rm dom}\,f={\rm dom}\,\phi\cap\mathcal{S}$. Since the collection of cells is finite, there exist an index $\ell\in[N]$ and an infinite index set $\mathcal{K}\subset\mathbb{N}$ such that $z^k\in C_\ell$ for all $k\in\mathcal{K}$. Clearly, $z\in{\rm cl}\,C_{\ell}$. Since $C_\ell\cap{\rm dom}\,\phi\neq\varnothing$, there exists some $i_0\in[m]$ such that $C_{\ell}\cap P_{i_0}\ne\varnothing$, which by  Corollary~\ref{cor:PLQ-consistency}(i) implies $C_\ell\subseteq P_{i_0}$. Consequently, $C_{\ell}\subset{\rm dom}\,\phi$, and $\ell\in\mathcal{J}$. Moreover, $z\in{\rm cl}\,C_\ell$ and $z\in C_j$. By Lemma~\ref{lem:hplane-partition}(i)-(ii), we have
 \begin{equation}\label{eq:jlrelation-PLQ}
 C_j\subset{\rm cl}\,C_{\ell},\ \mathbb{A}_j\subseteq \mathbb{A}_{\ell}\ \ {\rm and}\ \ \mathbb{L}_j\subseteq \mathbb{L}_{\ell}.
 \end{equation} 
Furthermore, by Corollary~\ref{cor:PLQ-consistency}(iii), 
$\langle (\widetilde{\mathcal{Q}}_j-\widetilde{\mathcal{Q}}_\ell)z+\widetilde{a}_j-\widetilde{a}_\ell,d\rangle=0$ for all $d\in\mathbb{L}_j$.
Then
\begin{equation}\label{eq:gradagree-PLQ}
\langle\widetilde{\mathcal{Q}}_jz+\widetilde{a}_j,d\rangle=\langle \widetilde{\mathcal{Q}}_{\ell}z+\widetilde{a}_{\ell},d\rangle\quad\  \forall d\in\mathbb{L}_j.
\end{equation}
To establish the desired equality, it suffices to prove that, for all $d\in T_z\mathcal{M}_j$, 
\begin{equation}\label{aim-equa}
\langle v,d\rangle=\langle\widetilde{\mathcal{Q}}_jz+\widetilde{a}_j,d\rangle.
\end{equation} 
Fix any $d\in T_z\mathcal{M}_j$. If $d=0$, \eqref{aim-equa} holds immediately. Suppose that $d\ne0$. By the preceding argument, $r_j>0$ and $T_z\mathcal{M}_j=\mathbb{L}_j\cap z^\perp$. 
Let $p_{\ell}$ be the orthogonal projection of the origin onto
$\mathbb{A}_{\ell}$. Then $\mathbb{A}_{\ell}=p_{\ell}+\mathbb L_\ell$ and $p_{\ell}\perp\mathbb{L}_{\ell}$. By \eqref{eq:jlrelation-PLQ}, we have $z,p_j\in\mathbb{A}_{\ell}$ and
$d\in\mathbb{L}_{\ell}\cap z^\perp$. 
It follows that $z\in\mathcal M_\ell$ and $\|p_\ell\|\le\|p_j\|<1$.
In particular, $r_\ell:=\sqrt{1-\|p_\ell\|^2}>0$.
For each $k\in\mathcal K$, let $d^k$ be the orthogonal projection
of $d$ onto $T_{z^k}\mathcal{M}_\ell=\mathbb{L}_\ell\cap(z^k-p_\ell)^\perp$. Recall that $z^k\in C_{\ell}\cap\mathcal{S}\subset\mathbb{A}_{\ell}\cap\mathcal{S}$ for every $k\in\mathcal{K}$. We have $\|z^k-p_\ell\|=r_\ell$, and hence
\[
d^k=d-\frac{\langle d,z^k-p_\ell\rangle}{r_\ell^2}(z^k-p_\ell)=d-\frac{\langle d,z^k-z\rangle}{r_\ell^2}(z^k-p_\ell)\quad\forall k\in\mathcal{K}.
\]
where the second equality follows from $d\perp p_\ell$ and $d\perp z$. Thus $\lim_{\mathcal{K}\ni k\to\infty}d^k=d$. For each $k\in\mathcal{K}$, applying \eqref{Fsubdiff-PLQ} at $z^k$ to  
$v^k\in\widehat\partial f(z^k)$ and taking the inner product with
$d^k\in T_{z^k}\mathcal{M}_\ell$, we obtain
$\langle v^k,d^k\rangle=\langle\widetilde{\mathcal Q}_\ell z^k+\widetilde{a}_\ell,d^k\rangle$. Passing to the limit $\mathcal{K}\ni k\to\infty$ in the equality yields
$\langle v,d\rangle =\langle\widetilde{\mathcal{Q}}_\ell z+\widetilde{a}_\ell,d\rangle$. Together with \eqref{eq:gradagree-PLQ}, we have $\langle v,d\rangle =\langle\widetilde{\mathcal{Q}}_jz+\widetilde{a}_j,d\rangle$. Since $d$ is arbitrary in $T_z\mathcal{M}_j$, equality \eqref{aim-equa} follows.
\end{proof}
\begin{remark}\label{remark-tangent}
The tangent-projection identity established in 
Proposition~\ref{prop:tangent-proj-PLQ} is closely related to the projection formula for stratifiable functions in \cite[Proposition~4]{Bolte2007-II}. Our proof uses the finite polyhedral cell decomposition associated with the PLQ representation of $\phi$ and the
geometry of the unit sphere to obtain an explicit formula
for the tangential projection of each limiting subgradient
of $f$.
\end{remark}

Li and Pong~\cite[Theorem~3.7]{LiPong2018} showed that a
tangential gradient inequality on the active manifold yields
a KL inequality with the same exponent for the original
function under prox-regularity, $C^2$ partial smoothness,
and nondegeneracy. Wu et al.~\cite[Theorem~5.2(a)]{WuLeeXie2026}
showed that a KL inequality in the ambient space with a
concave desingularizing function passes to the restriction
to the active manifold under $C^1$ partial smoothness and
nondegeneracy. In our setting, the following lemma uses the projection
identity in Proposition~\ref{prop:tangent-proj-PLQ} to
combine tangential gradient estimates with a common
exponent $\theta$ on $C_j\cap\mathcal{S}$ for the cells
$C_j$ whose closures contain the reference point. This
yields a KL inequality for $f$ with exponent $\theta$
at that point. 

\begin{lemma}\label{lem:cellwise-KL}
In the setting of Proposition~\ref{prop:tangent-proj-PLQ}, fix any
$\overline{z}\in{\rm dom}\,f$ and $\theta\!\in[0,1)$, and let
$\mathcal{I}(\overline{z}):=\{j\in\mathcal{J}\mid \overline{z}\in{\rm cl}\,C_j\}$. Suppose that for each $j\in\mathcal{I}(\overline{z})$ there exist a neighborhood $\mathcal{U}_j$ of $\overline{z}$ and constants $c_j>0$ and $\eta_j>0$ such that
\[
  \|\mathcal{P}_{T_z\mathcal{M}_j}\nabla q_j(z)\|\ge c_j(q_j(z)-q_j(\overline{z}))^\theta
\]
for all $z\in\mathcal{U}_j\cap C_j\cap\mathcal{S}$ with $0<q_j(z)-q_j(\overline{z})<\eta_j$. Then $f$ has the KL property with exponent $\theta$ at $\overline{z}$.
\end{lemma}
\begin{proof}
By Lemma~\ref{lem:domain-Lsubdiff}, $\overline{z}\in{\rm dom}\,\partial f$. If a cell $C_j$ meets ${\rm dom}\,\phi$, it meets some $P_i$.
By Corollary~\ref{cor:PLQ-consistency}(i), we then have
$C_j\subseteq P_i\subseteq{\rm dom}\,\phi$, so $j\in\mathcal{J}$.
Since $C_1,\ldots,C_N$ partition $\mathbb{E}$, the cells
indexed by $\mathcal{J}$ partition ${\rm dom}\,\phi$. 
Then $\mathcal{I}(\overline{z})$ is nonempty.
Define  
\[
 \mathcal{U}:=\Big(\bigcap_{j\in\mathcal{I}(\overline{z})}\mathcal{U}_j\Big) \setminus\Big(\bigcup_{j\in\mathcal{J}\setminus\mathcal{I}(\overline{z})}{\rm cl}\,C_j\Big),\ c:=\min_{j\in\mathcal{I}(\overline{z})}c_j\ \ {\rm and}\ \ \eta:=\min_{j\in\mathcal{I}(\overline{z})}\eta_j.
\]
Since $\mathcal{J}$ is finite and the excluded closures do not contain
$\overline{z}$, the above $\mathcal{U}$ is a neighborhood of $\overline{z}$ and $c,\eta>0$. Take any $z\in\mathcal{U}$ satisfying $0<f(z)-f(\overline{z})<\eta$.
Then $z\in{\rm dom}\,f={\rm dom}\,\phi\cap\mathcal{S}$. 
Since the cells $C_j$ for $j\in\mathcal{J}$ partition
${\rm dom}\,\phi$, there exists a unique $j_0\in\mathcal{J}$
such that $z\in C_{j_0}$. If $j_0\notin\mathcal{I}(\overline{z})$,
then $\mathcal{U}\cap{\rm cl}\,C_{j_0}=\varnothing$ by the definition of $\mathcal{U}$, contradicting $z\in\mathcal{U}\cap C_{j_0}$. Hence
$z\in C_{j_0}\cap\mathcal{S}$ for $j_0\in\mathcal{I}(\overline{z})$. 
Since $P_{i_{j_0}}$ is closed and contains $C_{j_0}$, we have $\phi=q_{j_0}$ on ${\rm cl}\,C_{j_0}$ and 
$0<q_{j_0}(z)-q_{j_0}(\overline{z})=f(z)-f(\overline{z})<\eta\le\eta_{j_0}$. Together with the hypothesis of the lemma and $c_{j_0}\ge c$, it follows that 
\[
\bigl\|\mathcal P_{T_z\mathcal{M}_{j_0}}\nabla q_{j_0}(z)\bigr\|
\ge c_{j_0}\bigl(f(z)-f(\overline z)\bigr)^\theta
\ge c\bigl(f(z)-f(\overline z)\bigr)^\theta.
\]
By Proposition~\ref{prop:tangent-proj-PLQ}, every $v\in\partial f(z)$ satisfies
$\mathcal{P}_{T_z\mathcal{M}_{j_0}}v =\mathcal{P}_{T_z\mathcal{M}_{j_0}}\nabla q_{j_0}(z)$. Then, 
\[
 \|v\|\ge \bigl\|\mathcal{P}_{T_z\mathcal{M}_{j_0}}v\bigr\|
	= \bigl\|\mathcal{P}_{T_z\mathcal{M}_{j_0}}\nabla q_{j_0}(z)\bigr\|\ge c\bigl(f(z)-f(\overline z)\bigr)^\theta.
\]
Taking the infimum over $v\in\partial f(z)$ yields
${\rm dist}(0,\partial f(z))\ge c\bigl(f(z)-f(\overline{z})\bigr)^\theta$. This proves that $f$ satisfies the KL property with exponent $\theta$ at $\overline{z}$.
\end{proof}
\subsection{KL exponent $3/4$ for PLQ functions}\label{sec3.2}
By Lemma~\ref{lem:cellwise-KL}, local gradient estimates for
the quadratic pieces of $\phi$ on the corresponding affine
sections of $\mathcal{S}$ yield a KL inequality for $f$ at a
reference point. The following theorem expresses an admissible
KL exponent for $f$ in terms of the smallest KL exponents on
these sections, as determined in Proposition~\ref{prop:quadratic-sphere-KL} below.
\begin{theorem}\label{theorem:PLQ-34KL}
Let $\phi$ and its PLQ representation be as in Corollary~\ref{cor:PLQ-consistency}, and suppose that ${\rm dom}\,\phi\cap\mathcal{S}\ne\varnothing$. 
Then the following statements hold for the function $f$. 
\begin{itemize}
\item[(i)] Use the notation of Proposition~\ref{prop:tangent-proj-PLQ} and Lemma~\ref{lem:cellwise-KL}. Fix any
$\overline{z}\in{\rm dom}\,f$. For each 
$j\in\mathcal{I}(\overline{z})$, let $\theta_j(\overline{z})$
denote the smallest KL exponent of $q_j+\delta_{\mathcal{M}_j}$
at $\overline{z}$. Then $f$ satisfies the KL property at
$\overline{z}$ with exponent $\theta_{\overline{z}}
:=\max_{j\in\mathcal{I}(\overline{z})}\theta_j(\overline{z})$.

\item[(ii)] For every $\overline{z}\in{\rm dom}\,f$, the
exponent in \textup{(i)} satisfies $\theta_{\overline{z}}\le3/4$.
Consequently, $f$ satisfies the KL property with exponent
$3/4$ throughout ${\rm dom}\,f$.
\end{itemize}
\end{theorem}

To determine the exponents $\theta_j(\overline{z})$ in item~(i)
of Theorem~\ref{theorem:PLQ-34KL}, we study the quadratic pieces
of $\phi$ restricted to their corresponding affine sections of
$\mathcal{S}$. Let $q$ be one of these pieces, and suppose that its associated section $\mathcal{M}$ has dimension $d\ge 1$. At a critical point
$\overline{z}$ of $q+\delta_{\mathcal{M}}$, stereographic
coordinates $z=\Psi(u)$ with $\Psi(0)=\overline{z}$ give
$q(\Psi(u))-q(\overline{z})=\rho(u)Q(u)$, where $\rho$ is
smooth and positive near $0$ and $Q$ is a polynomial of degree
at most four. The smallest KL exponent of $q+\delta_{\mathcal{M}}$
at $\overline{z}$ equals the smallest admissible exponent in
the one-sided gradient inequality for $Q$ at $0$. For polynomials of degree at most four in $d$ variables, the estimate in~\cite[Theorem~4.2]{DAcuntoKurdyka2005} provides an upper bound of $1-1/(4\cdot9^{d-1})$ on the
{\L}ojasiewicz gradient exponent at a critical zero. The following lemma uses the particular form of $Q$ to determine the smallest admissible exponent in the one-sided inequality and shows that it is at most $3/4$, independently of $d$.
\begin{lemma}\label{lem:quartic-gradient}
Let $\mathbb{V}$ be a nonzero finite-dimensional Euclidean space, and define $Q(u):=\langle u,\mathcal{H} u\rangle +\langle b,u\rangle\|u\|^2+\gamma\|u\|^4$ for $u\in\mathbb{V}$, where $\mathcal{H}:\mathbb{V}\to\mathbb{V}$ is a self-adjoint linear operator, $b\in\mathbb{V}$ and $\gamma\in\mathbb{R}$. Call $\theta\in[0,1)$ an admissible KL exponent of $Q$ at $0$ if there exist $\kappa>0$ and $\varepsilon>0$ such that 
$Q(u)^\theta\le\kappa\|\nabla Q(u)\|$ for all $u\in\mathbb{B}_{\mathbb{V}}^{\circ}(0,\varepsilon)$ with $Q(u)>0$. Let $\mathbb{M}:=\ker\mathcal{H}$, $\mathbb{L}\!:=\mathbb{M}^\perp$ and $b_{\mathbb{M}}\!:=\mathcal{P}_{\mathbb{M}}b$. The following statements hold. 
\begin{itemize}
\item[(i)] If $b_{\mathbb{M}}\ne0$, then the smallest admissible exponent is $2/3$.

\item[(ii)] Suppose that $b_{\mathbb{M}}=0$. Let $c_{\mathbb{L}}$ be the unique vector in $\mathbb{L}$ satisfying $2\mathcal{H} c_{\mathbb{L}}=b$, and set $\widehat{\gamma}:=\gamma-\langle c_{\mathbb{L}},\mathcal{H}c_{\mathbb{L}}\rangle$. If $\mathbb{M}\ne\{0\}$ and $\widehat{\gamma}>0$, then the smallest admissible exponent is $3/4$. In the remaining cases, it is $1/2$ if $\mathcal H$ has a positive eigenvalue and $0$ otherwise.
\end{itemize}
\end{lemma}
\begin{proof}
For each $u\in\mathbb{V}$, write $u=x+y$ with  
$x:=\mathcal{P}_{\mathbb{M}}u$ and $y:=\mathcal{P}_{\mathbb{L}}u$, and let $s:=\|u\|$ and $\chi:=\|\nabla Q(u)\|$. The self-adjointness of $\mathcal{H}$ implies that $\mathcal{H}(\mathbb{L})\subseteq\mathbb{L}$. Define $\mathcal{H}_{\mathbb{L}}:=\mathcal{H}|_{\mathbb{L}}
:\mathbb{L}\to\mathbb{L}$. If $\mathbb{L}\ne\{0\}$, $\mathcal{H}_{\mathbb{L}}$ is invertible. Then there exists $\sigma>0$ such that 
\begin{equation}\label{nonsingular-ineq}
\|\mathcal{H}y'\|\ge\sigma\|y'\|\quad\ \forall y'\in\mathbb{L}. 
\end{equation}
When $\mathbb{L}=\{0\}$, \eqref{nonsingular-ineq} holds with $\sigma=1$. Differentiating $Q$ and using $\mathcal{H}x=0$ gives
\[
 \nabla Q(u)=2\mathcal{H}y+bs^2+\bigl(2\langle b,u\rangle+4\gamma s^2\bigr)u\quad\ \forall u\in\mathbb{V}.
\]
Projecting this equality onto $\mathbb{L}$ and using the bound in \eqref{nonsingular-ineq} on $\mathcal{H}y$, we obtain
\[
2\sigma\|y\|\le \chi+\|b_{\mathbb{L}}\|s^2 +\bigl(2\|b\|s+4|\gamma|s^2\bigr)\|y\|\quad {\rm with}\ b_{\mathbb{L}}:=\mathcal{P}_{\mathbb{L}}b.
\]
For sufficiently small $u$, the coefficient of $\|y\|$ on the right-hand side is at most $\sigma$, so  
\begin{equation}\label{eq:quartic-y-bound}
\|y\|\le \alpha_0(\chi+s^2)\ \ {\rm with}\ \ \alpha_0:=\sigma^{-1}\max\{1,\|b_{\mathbb{L}}\|\}.
\end{equation}
In addition, since $\nabla Q$ is continuous and $\nabla Q(0)=0$, we may shrink the neighborhood of the origin so that $s\le1$ and
$\chi\le1$ throughout the following estimates. 

\smallskip
\noindent
\textbf{(i)} Taking the inner product of $\nabla Q(u)$ with $b_{\mathbb{M}}$, we obtain 
\[
 \langle b_{\mathbb{M}},\nabla Q(u)\rangle =\|b_{\mathbb{M}}\|^2s^2+2\langle b_{\mathbb{M}},x\rangle^2+2\langle b_{\mathbb{L}},y\rangle\langle b_{\mathbb{M}},x\rangle+4\gamma s^2\langle b_{\mathbb{M}},x\rangle\quad\ \forall u\in\mathbb{V}.
\]
Since the second term is nonnegative, the bounds
$\langle b_{\mathbb{M}},\nabla Q(u)\rangle
\le\|b_{\mathbb{M}}\|\chi$ and $\|x\|\le s$ give
\[
 \|b_{\mathbb{M}}\|^2s^2\le \|b_{\mathbb{M}}\|\chi+2\|b_{\mathbb{M}}\|\|b_{\mathbb{L}}\|\|y\|s+4|\gamma|\|b_{\mathbb{M}}\|s^3.
\]
Substituting \eqref{eq:quartic-y-bound} into this inequality and using $s\le1$ yields, for small enough $\|u\|$,
\[
 \|b_{\mathbb{M}}\|^2s^2\le \|b_{\mathbb{M}}\|(1+2\|b_{\mathbb{L}}\|\alpha_0)\chi+\|b_{\mathbb{M}}\|(2\|b_{\mathbb{L}}\|\alpha_0+4|\gamma|)s^3.
\]
Since $b_{\mathbb{M}}\ne0$, shrinking the neighborhood of the origin yields $(2\|b_{\mathbb{L}}\|\alpha_0+4|\gamma|)s\le\frac{\|b_{\mathbb{M}}\|}{2}$ for all $u$ in the neighborhood. Then $s^2\le\alpha_1\chi$ for $\alpha_1=\frac{2(1+2\|b_{\mathbb{L}}\|\alpha_0)}{\|b_{\mathbb{M}}\|}$. Substituting this into \eqref{eq:quartic-y-bound} yields $\|y\|\le \alpha_0(1+\alpha_1)\chi$ for all such $u$. Since $\langle u,\mathcal{H}u\rangle=\langle y,\mathcal{H}y\rangle$, we have
\[
|Q(u)|\le\|\mathcal{H}\|\|y\|^2+\|b\|s^3+|\gamma|s^4\quad\ \forall u\in\mathbb{V}.
\]
Note that $s^3\le \alpha_1^{3/2}\chi^{3/2}$, $s^4\le s^3$, and $\chi^2\le\chi^{3/2}$ for all sufficiently small $u$. The above inequality implies that $|Q(u)|\le \overline{\alpha}\chi^{3/2}$ with some $\overline{\alpha}>0$ for all sufficiently small $u$. Taking the power $2/3$ at points where $Q(u)>0$ proves that $2/3$ is admissible.

For minimality, take $v:=b_{\mathbb{M}}/\|b_{\mathbb{M}}\|$. Since $\mathcal{H}v=0$ and $\langle b,v\rangle=\|b_{\mathbb{M}}\|$, for every $t>0$, we have
$Q(tv)=\|b_{\mathbb{M}}\|t^3+\gamma t^4$ and $\nabla Q(tv)=t^2(b+2\|b_{\mathbb{M}}\|v)+4\gamma t^3v$. Moreover, $Q(tv)\ge(\|b_{\mathbb{M}}\|/2)t^3>0$ for all sufficiently small $t>0$. Hence, for each $\theta\in[0,2/3)$,
\[
\frac{\|\nabla Q(tv)\|}{Q(tv)^\theta}
=O(t^{2-3\theta})\to 0\quad\text{as}\ t\downarrow0.
\]
This limit rules out all $\theta\in[0,2/3)$, since admissibility
requires the ratio to be bounded below by a positive constant.
Consequently, the smallest admissible exponent is $2/3$.

\smallskip
\noindent
\textbf{(ii)} In this case, $b\in\mathbb{L}$. If $\mathbb{L}\ne\{0\}$, the invertibility of $\mathcal{H}_{\mathbb{L}}$ ensures that $c_{\mathbb{L}}$ is the unique vector in $\mathbb{L}$ such that $2\mathcal{H}c_{\mathbb{L}}=b$. If $\mathbb{L}=\{0\}$, then $b=0$ and the unique solution in $\mathbb{L}$ is $c_{\mathbb{L}}=0$. Set $w:=y+c_{\mathbb{L}}s^2$. Expanding the quadratic form $\langle w,\mathcal{H} w\rangle$ then gives
\[
\langle w,\mathcal{H} w\rangle
 =\langle y,\mathcal{H}y\rangle +2s^2\langle y,\mathcal{H}c_{\mathbb{L}}\rangle+s^4\langle c_{\mathbb{L}},\mathcal{H}c_{\mathbb{L}}\rangle=\langle y,\mathcal{H} y\rangle +s^2\langle b,y\rangle+s^4\langle c_{\mathbb{L}},\mathcal {H}c_{\mathbb{L}}\rangle.
\]
Since $\langle u,\mathcal{H}u\rangle=\langle y,\mathcal{H}y\rangle$ and $\langle b,u\rangle=\langle b,y\rangle$, the above equality implies that
\begin{equation}\label{eq:quartic-square-completion}
 Q(u)=\langle w,\mathcal{H}w\rangle+\widehat{\gamma}s^4\quad\ \forall u\in\mathbb{V}.
\end{equation}
Moreover, noting that $2\mathcal{H}y+bs^2=2\mathcal{H}w$ and $\langle b,w\rangle=\langle b,y\rangle+2\langle c_{\mathbb{L}},\mathcal{H}c_{\mathbb{L}}\rangle s^2$, we have $2\langle b,u\rangle+4\gamma s^2=2\langle b,y\rangle+4\gamma s^2=2\langle b,w\rangle+4\widehat{\gamma}s^2$. Therefore, it holds
\[
 \nabla Q(u)=2\mathcal{H}w+(2\langle b,w\rangle+4\widehat{\gamma} s^2) u\quad\ \forall u\in\mathbb{V}.
\]
As $\mathcal{H}w\in\mathbb L$, projecting the above equality onto $\mathbb{L}$ and $\mathbb{M}$, respectively, leads to 
\begin{equation}\label{eq:quartic-gradient-projections}
\mathcal{P}_{\mathbb{L}}\nabla Q(u)=2\mathcal{H} w+(2\langle b,w\rangle+4\widehat{\gamma} s^2)y\ \ {\rm and}\ \ 
\mathcal{P}_{\mathbb{M}}\nabla Q(u)=(2\langle b,w\rangle+4\widehat{\gamma} s^2)x.
\end{equation}

\noindent
\textbf{Case 1: $\mathbb{M}=\{0\}$ or $\widehat{\gamma}\le0$.}
First suppose that $\mathcal{H}$ has no positive eigenvalue.
If $\widehat{\gamma}\le0$, then \eqref{eq:quartic-square-completion}
gives $Q(u)\le0$. If $\mathbb M=\{0\}$, then in this case $\mathcal{H}$ is negative definite, so there exists $\tau>0$ such that $Q(u)\le -\tau s^2+\|b\|s^3+|\gamma|s^4\le -({\tau}/{2})s^2$  for all small enough $u$. In either case, $Q(u)\le0$ for all sufficiently small $u$.
Thus the defining inequality holds vacuously with exponent $0$,
the smallest admissible one. Now suppose that $\mathcal{H}$ has a positive eigenvalue. We argue that the exponent $1/2$ is admissible.

\noindent
\textbf{Subcase 1.1: $\mathbb{M}=\{0\}$.} In this case, $y=u$. For all sufficiently small $u$, we have $\alpha_0s\le 1/2$ and
$s\le1$. Combining $\alpha_0s\le 1/2$ with \eqref{eq:quartic-y-bound} leads to $s\le2\alpha_0\chi$. Set $\alpha:=\|\mathcal H\|+\|b\|+|\gamma|$. For all such $u$ with $Q(u)>0$, the bound $|Q(u)|\le \alpha s^2$ gives
\[
Q(u)^{1/2}\le \sqrt{\alpha}\,s\le 2\alpha_0\sqrt{\alpha}\,\chi
=2\alpha_0\sqrt{\alpha}\,\|\nabla Q(u)\|.
\]
Consequently, the exponent $1/2$ is admissible.  

\noindent
\textbf{Subcase 1.2: $\mathbb{M}\ne\{0\}$.} Now $\widehat{\gamma}\le0$. At any point with $Q(u)>0$, it follows from \eqref{eq:quartic-square-completion} that $|\widehat{\gamma}|s^4\le \langle w,\mathcal{H}w\rangle
\le \|\mathcal{H}\|\|w\|^2$, and hence
$|\widehat{\gamma}|s^2\le \sqrt{|\widehat{\gamma}|\|\mathcal{H}\|}\|w\|$.
Consequently,
\[
 |2\langle b,w\rangle+4\widehat{\gamma} s^2|\le \gamma_0\|w\|\ \ {\rm with}\ \ 
\gamma_0:=2\|b\|+4\sqrt{|\widehat\gamma|\|\mathcal H\|}.
\]
The first identity in \eqref{eq:quartic-gradient-projections}, together with \eqref{nonsingular-ineq}, then gives $2\sigma\|w\|\le \chi+\gamma_0\|y\|\|w\|$. Since $y\to0$ as $u\to0$, we may take $\gamma_0\|y\|\le\sigma$. Thus $\|w\|\le \chi/\sigma$. Using $\widehat\gamma\le0$ in \eqref{eq:quartic-square-completion} yields
\[
0<Q(u)\le\|\mathcal{H}\|\|w\|^2\le({\|\mathcal{H}\|}/{\sigma^2})\chi^2
\]
for all sufficiently small $u$ with $Q(u)>0$. This again gives admissibility of $1/2$. 

To prove minimality, let $v$ be a unit eigenvector of $\mathcal{H}$
associated with a positive eigenvalue $\lambda$.
As $t\downarrow0$, we have $Q(tv)=\lambda t^2+O(t^3)$ and $\nabla Q(tv)=2\lambda tv+O(t^2)$. Then $Q(tv)\ge(\lambda/2)t^2>0$ for all sufficiently small $t>0$. However, for every $\theta\in[0,1/2)$,
\[
 \|\nabla Q(tv)\|/Q(tv)^\theta=O(t^{1-2\theta})\to0\ \ {\rm as}\ t\downarrow 0.
\]
This shows that $1/2$ is the smallest admissible exponent.

\noindent
\textbf{Case 2: $\mathbb M\ne\{0\}$ and $\widehat\gamma>0$.}
From \eqref{nonsingular-ineq} and the first identity in \eqref{eq:quartic-gradient-projections}, we have
\[
2\sigma\|w\|\le\chi+\bigl(2\|b\|\|w\|+4\widehat{\gamma} s^2\bigr)\|y\|.
\]
For $u$ sufficiently close to zero, $2\|b\|\|y\|\le\sigma$.
Combining the resulting estimate with \eqref{eq:quartic-y-bound}
and $s\le1$ gives $\sigma\|w\|\le\chi+4\widehat{\gamma} s^2\|y\|
\le (1+4\widehat{\gamma}\alpha_0)\chi+4\widehat{\gamma}\alpha_0s^4$. In particular, 
\begin{equation}\label{temp-wineq}
 \|w\|\le \alpha_2(\chi+s^4)\ \ {\rm with}\ \alpha_2:=\sigma^{-1}(1+4\widehat{\gamma}\alpha_0).
\end{equation}
The second identity in \eqref{eq:quartic-gradient-projections}
gives $|2\langle b,w\rangle+4\widehat{\gamma} s^2|\|x\|\le\chi$, which in turn implies that $4\widehat{\gamma}s^2\|x\|\le\chi+2\|b\|\|w\|\|x\|\le \chi+2\|b\|\|w\|s$. Together with $s\le\|x\|+\|y\|$, we obtain
\begin{align*}
4\widehat{\gamma}s^3\le\chi+2\|b\|\|w\|s+4\widehat{\gamma} s^2\|y\|
&\overset{\eqref{temp-wineq},\eqref{eq:quartic-y-bound}}{\le} \chi+2\|b\|\alpha_2s(\chi+s^4)+4\widehat{\gamma}\alpha_0s^2(\chi+s^2)\\
&\le (1+\alpha_3)\chi+\alpha_3s^4\ \ {\rm with}\ \alpha_3:=2\|b\|\alpha_2+4\widehat{\gamma}\alpha_0,
\end{align*}
where the last inequality uses $s\le1$. Since $\widehat{\gamma}>0$, $\alpha_3s\le2\widehat{\gamma}$ for all sufficiently small $u$, so we have 
$2\widehat{\gamma}s^3\le(1+\alpha_3)\chi$. This gives $s^3=O(\chi)$ and, since $s^4\le s^3$, $\|w\|=O(\chi)$ follows from \eqref{temp-wineq}. Together with the previous \eqref{eq:quartic-square-completion} and $\chi\le1$, we obtain 
\[
 |Q(u)|\le\|\mathcal{H}\|\|w\|^2+\widehat{\gamma}s^4       =O(\chi^2+\chi^{4/3})=O(\chi^{4/3}). 
\]
Taking the power $3/4$ at points with $Q(u)>0$ proves admissibility of $3/4$.

Finally, fix any unit vector $v\in\mathbb{M}$ and set
$u(t):=tv-c_{\mathbb{L}}t^2$ for $t\in\mathbb{R}$. Note that $v\perp c_{\mathbb{L}}$. Along this curve, we have $y=\mathcal{P}_{\mathbb{L}}u(t)=-c_{\mathbb{L}}t^2$, $s^2=\|u(t)\|^2=t^2+\|c_{\mathbb{L}}\|^2t^4$ and $w=y+c_{\mathbb{L}}s^2=\|c_{\mathbb{L}}\|^2t^4c_{\mathbb{L}}$. 
By \eqref{eq:quartic-square-completion},
$Q(u(t))=\widehat\gamma t^4+O(t^6)\ge (\widehat{\gamma}/{2})t^4>0$ for sufficiently small $t>0$. Moreover, $\|w\|=O(t^4)$,
$2\langle b,w\rangle+4\widehat{\gamma} s^2=O(t^2)$, and $\|u(t)\|=O(t)$. The identity
$\nabla Q(u(t))=2\mathcal H w+(2\langle b,w\rangle+4\widehat{\gamma} s^2) u(t)$ thus gives $\|\nabla Q(u(t))\|=O(t^3)$. Then, for every $\theta\in[0,3/4)$, ${\|\nabla Q(u(t))\|}/{Q(u(t))^\theta}
=O(t^{3-4\theta})\to 0$ as $t\downarrow0$. This shows that $3/4$ is the smallest admissible exponent. We complete the proof.
\end{proof}

We next use Lemma~\ref{lem:quartic-gradient} to determine the
smallest KL exponent of $q+\delta_{\mathcal{M}}$ at each point
of $\mathcal{M}$, where $q$ is quadratic and $\mathcal{M}$ is an
affine section of $\mathcal{S}$.
\begin{proposition}\label{prop:quadratic-sphere-KL}
Let $\mathcal{Q}:\mathbb{E}\to\mathbb{E}$ be a self-adjoint linear operator, and define
\[
 q(z):=\tfrac12\langle z,\mathcal{Q}z\rangle+\langle a,z\rangle+\beta\quad{\rm for}\ z\in\mathbb{E}, 
\]
where $a\in\mathbb{E}$ and $\beta\in\mathbb{R}$. Let $\mathbb{L}\subseteq\mathbb{E}$ be a linear subspace and $p\in\mathbb{L}^\perp$ such that $\mathcal{M}:=(p+\mathbb{L})\cap\mathcal{S}$ is nonempty. Set $h:=q+\delta_{\mathcal{M}}$. Fix any $\overline{z}\in\mathcal{M}$ and let $r:=\|\overline{z}-p\|$. 
\begin{itemize}
\item[(i)] If $\overline{z}$ is isolated in $\mathcal{M}$ or
$0\notin\partial h(\overline{z})$, the smallest KL exponent
of $h$ at $\overline{z}$ is $0$.

\item[(ii)] Suppose that $\overline{z}$ is nonisolated in $\mathcal{M}$ and $0\in\partial h(\overline{z})$. Set $\overline{u}:=(\overline{z}-p)/r$ and $\mathbb{V}:=T_{\overline{z}}\mathcal{M}$. Define $\mu:={\langle\mathcal{Q}\overline{z}+a,\overline{u}\rangle}/r$ and $\mathcal{B}:=\mathcal{Q}-\mu\mathcal{I}$. Then, the smallest KL exponent of $h$ at $\overline{z}$ equals
the smallest admissible exponent in Lemma~\ref{lem:quartic-gradient}
with  $\mathcal{H}=\mathcal{P}_{\mathbb{V}}\mathcal{B}|_{\mathbb{V}},\,b=-2\mathcal{P}_{\mathbb{V}}\mathcal{B}\overline{u}$ and $\gamma=\langle\overline{u},\mathcal{B}\overline{u}\rangle$.
\end{itemize}
In particular, $3/4$ is an admissible KL exponent of $h$ at every point of $\mathcal{M}$.
\end{proposition}
\begin{proof}
Since $p\perp\mathbb{L}$, we have $\|p+\xi\|^2=\|p\|^2+\|\xi\|^2$ for every $\xi\in\mathbb{L}$. In particular, taking $\xi=\overline{z}-p$ and using $\|\overline{z}\|=1$ gives $\|p\|^2+r^2=1$. Thus, for every $\xi\in\mathbb{L}$, $\|p+\xi\|=1$ if and only if $\|\xi\|=r$. It follows that $\mathcal{M}=(p+\mathbb{L})\cap\mathcal{S}=\{p+\xi\mid \xi\in\mathbb{L},\,\|\xi\|=r\}$. If $r=0$, then $\mathcal{M}=\{p\}$. If $r>0$ and
$\dim\mathbb{L}=1$, then $\mathcal{M}$ consists of two isolated
points. In either case, $\overline{z}$ is isolated in
${\rm dom}\,h=\mathcal{M}$, so $\widehat{\partial}h(\overline{z})=\mathbb{E}$. Choose a neighborhood $\mathcal{U}_{\overline{z}}$ of $\overline{z}$ such that
$\mathcal{U}_{\overline{z}}\cap{\rm dom}\,h=\{\overline{z}\}$. Thus the KL inequality holds vacuously with exponent $0$, which is the smallest
admissible exponent. We therefore assume that $r>0$ and $\dim\mathbb{L}\ge2$. Now $\mathcal{M}$ is a smooth embedded manifold of dimension $\dim\mathbb{L}-1$. For $z\in\mathbb{E}$, set
$\widehat{u}(z):=(z-p)/r$. Then, for every $z\in\mathcal{M}$, 
$\widehat{u}(z)$ is a unit vector in $\mathbb{L}$, and
\begin{equation}\label{tangent-formula}
 T_z\mathcal{M}=\mathbb{L}\cap z^\perp
 =\mathbb{L}\cap\widehat{u}(z)^\perp\ \ {\rm and}\ \ 
 \mathcal{P}_{T_z\mathcal{M}} =\mathcal{P}_{\mathbb{L}}-\widehat{u}(z)\widehat{u}(z)^\top.
\end{equation}
Since $\mathcal{N}_{\mathcal{M}}(z)=(T_z\mathcal{M})^\perp$, we have
$\partial h(z)=\nabla q(z)+(T_z\mathcal{M})^\perp$ for all $z\in\mathcal{M}$. Hence 
\begin{equation}\label{subdiff-dist}
{\rm dist}\bigl(0,\partial h(z)\bigr)=\|G(z)\|\ \ {\rm with}\ \ G(z):=\mathcal{P}_{T_z\mathcal{M}}\nabla q(z)\ \ {\rm for}\ z\in\mathcal{M}.   
\end{equation}
Since $\widehat{u}(z)=(z-p)/r$, the projection formula in \eqref{tangent-formula} and the continuity of $\nabla q$ imply that $G$ is continuous on $\mathcal M$. We distinguish two cases.
	
\smallskip
\noindent
\textbf{Case 1: $G(\overline{z})\neq0$.} By the continuity of $G$, for all $z\in\mathcal{M}$ sufficiently close to $\overline{z}$, we have  ${\rm dist}(0,\partial h(z))=\|G(z)\|\geq\frac{1}{2}\|G(\overline{z})\|>0$. At such points with $h(z)>h(\overline{z})$, since $(h(z)-h(\overline{z}))^0=1$, the KL inequality holds
with exponent $0$, the smallest admissible one. 

\smallskip
\noindent
\textbf{Case 2: $G(\overline{z})=0$.} We divide the proof into the following three steps.

\noindent
\emph{Step 1. Stereographic parametrization and its local inverse.} With $\overline{u}$ and $\mathbb{V}$ as in the statement, taking $z=\overline{z}$ in \eqref{tangent-formula} gives $\mathbb{V}=\mathbb{L}\cap\overline{u}^{\perp}$ and $\mathcal{P}_{\mathbb{V}}=\mathcal{P}_{\mathbb{L}}-\overline{u}\overline{u}^\top$. Define 
\begin{equation}\label{Psi-fun}
\Psi(u):=p+r\frac{(1-\|u\|^2)\overline u+2u}{1+\|u\|^2}\ \ {\rm and}\ \ s:=\|u\|^2\quad{\rm for}\  u\in\mathbb{V}.
\end{equation} 
For any $u\in\mathbb{V}$, clearly, $\Psi(u)-p\in\mathbb{L}$. Moreover, since $u\perp\overline{u}$ and $\|\overline{u}\|=1$, we have
\[
\|\Psi(u)-p\|^2=r^2\frac{(1-\|u\|^2)^2+4\|u\|^2}{(1+\|u\|^2)^2}=r^2.
\]
By the representation of $\mathcal{M}$ above, it follows that
$\Psi(u)\in\mathcal{M}$. In particular, $\Psi(0)=p+r\overline{u}=\overline{z}$.
Set 
$t(z):=\langle\widehat{u}(z),\overline{u}\rangle$ for $z\in\mathbb{E}$.
Since $t(\overline{z})=1$, there exists an open neighborhood
$\mathcal{W}$ of $\overline{z}$ in $\mathbb{E}$ on which $1+t(z)>0$.
Define
\[
\mathcal{A}(z):=\frac{\mathcal{P}_{\mathbb{V}}\widehat{u}(z)}{1+t(z)}\quad{\rm for}\ z\in\mathcal{W}.
\]
Then $\mathcal{A}:\mathcal{W}\to\mathbb{V}$ is smooth and
$\mathcal{A}(\overline z)=0$. For any $z\in\mathcal{W}\cap\mathcal{M}$, we have
$\widehat{u}(z)\in\mathbb{L}$ and $\|\widehat{u}(z)\|=1$. Together with $\mathcal{P}_{\mathbb{V}}=\mathcal{P}_{\mathbb{L}}
-\overline{u}\,\overline{u}^\top$, we obtain
\[
\mathcal{A}(z)=\frac{\widehat{u}(z)-t(z)\overline{u}}{1+t(z)}\ \ {\rm and}\ \ \|\mathcal{A}(z)\|^2=\frac{1-t(z)}{1+t(z)}
\]
It follows that $\frac{1-\|\mathcal{A}(z)\|^2}{1+\|\mathcal{A}(z)\|^2}=t(z)$ and $\frac{2}{1+\|\mathcal{A}(z)\|^2}=1+t(z)$. Substituting these identities into \eqref{Psi-fun} yields $\Psi(\mathcal{A}(z))=p+r\bigl(t(z)\overline u+(1+t(z))\mathcal A(z)\bigr)=z$. Conversely, since $\Psi$ is continuous and $\Psi(0)=\overline{z}\in\mathcal{W}$, we have
$\Psi(u)\in\mathcal{W}$ for all $u$ sufficiently close to $0$.
For such $u$, substituting $\Psi(u)$ into the definition of $\widehat{u}(\cdot)$ yields 
\[
 \widehat{u}(\Psi(u))=\frac{\Psi(u)-p}{r}=\frac{(1-s)\overline u+2u}{1+s}.
\]
Since $\mathcal{P}_{\mathbb{V}}u=u$ and $\mathcal{P}_{\mathbb{V}}\overline{u}=0$, it follows that
$\mathcal{P}_{\mathbb{V}}\widehat{u}(\Psi(u))=\frac{2u}{1+s}$.
Moreover, $1+t(\Psi(u))=1+\langle\widehat u(\Psi(u)),\overline u\rangle
=\frac{2}{1+s}$. Hence $\mathcal A(\Psi(u))=u$ for all $u$ sufficiently close to $0$.

\smallskip
\noindent
\emph{Step 2. Reduction to a quartic polynomial.}
Set $\overline{g}:=\mathcal{Q}\overline{z}+a$. By the definition of $G$,  $G(\overline{z})=\mathcal{P}_{\mathbb{V}}\overline{g}$. Since
$G(\overline{z})=0$, we have $\mathcal{P}_{\mathbb{V}}\overline{g}=0$. The definition of $\mu$ gives $\langle\overline g,\overline u\rangle=r\mu$.
For $u\in\mathbb V$, set $w:=u-s\overline u$. By \eqref{Psi-fun}, we have $d:=\Psi(u)-\overline z=\frac{2r}{1+s}w$, $\langle\overline{g},w\rangle=-r\mu s$, and $\|w\|^2=s(1+s)$. Set $\rho(u)\!:=\frac{2r^2}{(1+s)^2}$ for $u\in\mathbb{V}$. The quadratic expansion of $q$ gives
\[
\widetilde{Q}(u):=q(\Psi(u))-q(\overline{z})
 =\langle\overline{g},d\rangle+\tfrac12\langle d,\mathcal{Q}d\rangle
 =\rho(u)\langle w,\mathcal{B}w\rangle\ \ {\rm for}\ u\in\mathbb{V}.
\]
Thus, $\widetilde{Q}(u)=\rho(u)Q(u)$, where
$Q(u):=\langle w,\mathcal{B}w\rangle=\langle u,\mathcal{H}u\rangle+\langle b,u\rangle\|u\|^2+\gamma\|u\|^4$, 
with $\mathcal{H},\,b$ and $\gamma$ as in the statement. 
Fix any $\theta\in[0,1)$. Since $\rho(u)>0$, the conditions
$\widetilde{Q}(u)>0$ and $Q(u)>0$ are equivalent. At such points,
using $\nabla\rho(u)=-\frac{4\rho(u)u}{1+s}$ gives
\[
\frac{\nabla\widetilde{Q}(u)}{\widetilde{Q}(u)^\theta}
=\rho(u)^{1-\theta}
\left(
\frac{\nabla Q(u)}{Q(u)^\theta}
-\frac{4Q(u)^{1-\theta}}{1+s}u
\right).
\]
For $u$ near $0$, the factor $\rho(u)^{1-\theta}$ is bounded above and away from zero. Since $Q(0)=0$ and $\theta<1$, the second term in parentheses tends to zero as $u\to0$ through points with
$Q(u)>0$. The triangle inequality in both directions shows that
the norm of the difference in parentheses has a positive lower
bound near $0$ on these points if and only if
$\|\nabla Q(u)\|/Q(u)^\theta$ does. Thus $\widetilde{Q}$ and $Q$
have the same admissible exponents at $0$.

\smallskip
\noindent\emph{Step 3. Transfer of the exponents to $h$.} For any $u,\xi\in\mathbb{V}$, the smooth curve $t\mapsto\Psi(u+t\xi)$ lies in $\mathcal{M}$ and passes through $\Psi(u)$ at $t=0$. Its velocity at $t=0$ is 
\[
\left.\frac{d}{dt}\Psi(u+t\xi)\right|_{t=0}
=\Psi'(u)\xi\in T_{\Psi(u)}\mathcal{M}.
\]
Since $\xi$ is arbitrary, we have $\Psi'(u)\mathbb{V}\subseteq T_{\Psi(u)}\mathcal{M}$. Using the chain rule then yields 
\[
\nabla\widetilde{Q}(u)=[\Psi'(u)]^*G(\Psi(u))\quad\ \forall u\in\mathbb{V}.
\]
Since $\Psi'$ is continuous at $0$, there exists $c_1>0$ such that, for all $u\in\mathbb{V}$ near $0$, 
\begin{equation}\label{temp-gradwtq}
\|\nabla\widetilde{Q}(u)\|\le c_1\|G(\Psi(u))\|.
\end{equation}
For $z\in\mathcal M$ near $\overline z$, the identity
$\Psi(\mathcal A(z))=z$ gives
$q(z)-q(\overline z)=\widetilde{Q}(\mathcal A(z))$.
Differentiating along a smooth curve in $\mathcal M$ through $z$
with velocity $\xi\in T_z\mathcal M$ yields
\[
\langle G(z),\xi\rangle
=\bigl\langle\nabla\widetilde{Q}(\mathcal{A}(z)),\mathcal{A}'(z)\xi\bigr\rangle.
\]
Since $\mathcal{A}'$ is bounded near $\overline{z}$, taking the
supremum over $\xi\in T_z\mathcal{M}$ with $\|\xi\|\le1$ gives
$\|G(z)\|\le c_2\|\nabla\widetilde{Q}(\mathcal{A}(z))\|$ for some $c_2>0$. Together with \eqref{temp-gradwtq}, for any $u\in\mathbb{V}$ sufficiently close to $0$, using $z=\Psi(u)$ and $\mathcal{A}(\Psi(u))=u$, we obtain
\begin{equation}\label{wtq-ineq}
c_2^{-1}\|G(\Psi(u))\|\le \|\nabla\widetilde{Q}(u)\|\le c_1\|G(\Psi(u))\|.
\end{equation}

Fix any $\theta\in[0,1)$. Suppose first that $\theta$ is admissible
for $\widetilde{Q}$ at $0$, and let $\kappa_1>0$ be a constant
in the corresponding gradient inequality. For $z\in\mathcal{M}$
sufficiently close to $\overline{z}$ with $h(z)>h(\overline{z})$,
set $u:=\mathcal{A}(z)$. Then $\Psi(u)=z$, so
$\widetilde{Q}(u)=h(z)-h(\overline{z})>0$. Since $\mathcal{A}(z)\to0$
as $z\to\overline{z}$, the gradient inequality for $\widetilde{Q}$
and the upper bound in \eqref{wtq-ineq} give
\[
\bigl(h(z)-h(\overline{z})\bigr)^\theta
=\widetilde{Q}(u)^\theta
\le \kappa_1\|\nabla\widetilde{Q}(u)\|
\le \kappa_1c_1\|G(z)\|
=\kappa_1c_1{\rm dist}(0,\partial h(z)).
\]
Hence $\theta$ is a KL exponent of $h$ at $\overline{z}$. Conversely, suppose that $\theta$ is a KL exponent of $h$ at $\overline{z}$. Let $\mathcal{U}_{\overline{z}}$ be the corresponding neighborhood,
$\eta>0$ the upper bound on the positive function value gap,
and $\kappa_2>0$ the constant multiplying the subdifferential
distance in the KL inequality. Since $\Psi(u)\to\overline{z}$
and $\widetilde{Q}(u)\to0$ as $u\to0$, for all sufficiently small
$u\in\mathbb{V}$ with $\widetilde{Q}(u)>0$, the point $z:=\Psi(u)$
belongs to $\mathcal{U}_{\overline{z}}$ and satisfies
$0<h(z)-h(\overline{z})=\widetilde{Q}(u)<\eta$.
The KL inequality for $h$ and the lower bound in
\eqref{wtq-ineq} therefore yield
\[
\widetilde{Q}(u)^\theta
\le \kappa_2{\rm dist}(0,\partial h(z))
=\kappa_2\|G(\Psi(u))\|
\le \kappa_2c_2\|\nabla\widetilde{Q}(u)\|.
\]
Thus $\theta$ is admissible for $\widetilde{Q}$ at $0$. Therefore, the smallest KL exponent of $h$ at $\overline{z}$ equals the smallest admissible exponent of $Q$ at $0$. Since $\dim\mathbb{V}=\dim\mathbb{L}-1\ge1$ and $\mathcal{H}$
is self-adjoint on $\mathbb{V}$, Lemma~\ref{lem:quartic-gradient}
gives the asserted classification. In all cases the smallest
exponent is at most $3/4$. Restricting the positive gap to be
less than $1$ therefore gives admissibility of $3/4$ at every
point of $\mathcal M$. Thus, the proof is completed.
\end{proof}

To prove Theorem~\ref{theorem:PLQ-34KL}(i), we use the local
tangential gradient estimates for the quadratic pieces $q_j$
of $\phi$ on $\mathcal{M}_j$ given by
Proposition~\ref{prop:quadratic-sphere-KL}. At each
$\overline{z}\in{\rm dom}\,f$, the maximum of the corresponding
KL exponents over $j\in\mathcal{I}(\overline{z})$ gives the
common exponent required by Lemma~\ref{lem:cellwise-KL}. Proposition~\ref{prop:quadratic-sphere-KL} bounds each section
exponent by $3/4$, so the pointwise exponent
$\theta_{\overline z}$ is at most $3/4$ at every point of
${\rm dom}\,f$.

\medskip
\begin{proof} [Proof of Theorem~\ref{theorem:PLQ-34KL}]
Fix any $\overline{z}\in{\rm dom}\,f={\rm dom}\,\partial f$, where the equality follows from Lemma \ref{lem:domain-Lsubdiff}. Since the cells indexed by $\mathcal{J}$ form a finite partition
of ${\rm dom}\,\phi$ and $\overline{z}\in{\rm dom}\,\phi$,
one of these cells contains $\overline{z}$. Hence
$\mathcal{I}(\overline{z})$ is finite and nonempty. For each $j\in\mathcal{I}(\overline{z})$, we have
$\overline{z}\in{\rm cl}\,C_j\cap\mathcal{S}
\subseteq\mathcal{M}_j$.
If $\overline{z}$ is isolated in $\mathcal{M}_j$, choose
a neighborhood $\mathcal{U}_j$ of $\overline{z}$ such that
$\mathcal{U}_j\cap\mathcal{M}_j=\{\overline{z}\}$. In this case,
the estimate below holds vacuously with
$c_j=\eta_j=1$. Otherwise,
Proposition~\ref{prop:quadratic-sphere-KL} gives a KL inequality
for $q_j+\delta_{\mathcal{M}_j}$ at $\overline{z}$ with exponent
$\theta_j(\overline{z})$. On $\mathcal M_j$, the function
$q_j+\delta_{\mathcal M_j}$ agrees with $q_j$, and its
subdifferential distance at $z$ equals
$\|\mathcal P_{T_z\mathcal M_j}\nabla q_j(z)\|$ by
\eqref{subdiff-dist}. Thus, in either case, there exist a neighborhood $\mathcal{U}_j$ of $\overline{z}$ and constants $c_j>0$ and $\eta_j\in(0,1]$ such that
\[
\bigl\|\mathcal{P}_{T_z\mathcal{M}_j}\nabla q_j(z)\bigr\|
\ge c_j\bigl(q_j(z)-q_j(\overline{z})\bigr)^{\theta_j(\overline{z})}
\]
for all $z\in\mathcal{U}_j\cap\mathcal{M}_j$ satisfying
$0<q_j(z)-q_j(\overline{z})<\eta_j$. Proposition~\ref{prop:quadratic-sphere-KL} gives $\theta_j(\overline{z})\le3/4$ for each
$j\in\mathcal{I}(\overline{z})$, so $\theta_{\overline{z}}\in[0,3/4]$.
Since $\theta_j(\overline{z})\le\theta_{\overline{z}}$ and the
positive gap is less than $1$, the same estimate holds with
$\theta_{\overline{z}}$ in place of $\theta_j(\overline{z})$.
Restricting these estimates to $C_j\cap\mathcal{S}$ for each
$j\in\mathcal{I}(\overline{z})$ verifies the hypotheses of
Lemma~\ref{lem:cellwise-KL} with $\theta=\theta_{\overline{z}}$.
The lemma gives the KL property of $f$ at $\overline{z}$ with
exponent $\theta_{\overline{z}}$, establishing \textup{(i)}.

Let $\mathcal U$ and $c,\eta>0$ be the neighborhood and constants
in the resulting KL inequality. Since $\theta_{\overline{z}}\le3/4$,
for all $z\in\!\mathcal U$ with $0<f(z)-\!f(\overline{z})<\min\{\eta,1\}$, we have
\[
{\rm dist}(0,\partial f(z))
\ge c\bigl(f(z)-f(\overline{z})\bigr)^{\theta_{\overline{z}}}
\ge c\bigl(f(z)-f(\overline{z})\bigr)^{3/4}.
\]
Since $\overline{z}$ is arbitrary in ${\rm dom}\,f$, this proves item \textup{(ii)}.
\end{proof}

The following example shows that the smallest KL exponent of
a quadratic function restricted to the unit circle can be
$2/3$ at a nonoptimal critical point.  
\begin{example}\label{exa:quadratic-sphere-two-thirds}
Let $\mathcal{S}$ be the unit circle in $\mathbb{R}^2$ and set
$h:=q+\delta_{\mathcal{S}}$, where $q(z):=z_2(1-z_1)$ for
$z=(z_1,z_2)^\top\in\mathbb{R}^2$. Define
$z(t):=(\cos t,\sin t)^\top$ for $t\in\mathbb{R}$ and set
$\overline{z}:=z(0)=(1,0)^\top$. Since $h=q$ on $\mathcal{S}$,
for every $t\in\mathbb{R}$ we have
\[
h(z(t))-h(\overline{z}) =q(z(t))-q(\overline{z})
=\sin t(1-\cos t).
\]
For every $t\in\mathbb{R}$, since $T_{z(t)}\mathcal{S}={\rm span}\{z'(t)\}$ and $\|z'(t)\|=1$, by \eqref{subdiff-dist} we have
\begin{align*}
{\rm dist}(0,\partial h(z(t)))
=\bigl\|\mathcal{P}_{T_{z(t)}\mathcal{S}}\nabla q(z(t))\bigr\|
=\bigl|\langle\nabla q(z(t)),z'(t)\rangle\bigr|
=\bigl|(1-\cos t)(1+2\cos t)\bigr|.
\end{align*}
The distance vanishes at $t=0$, and
$h(z(-t))<h(\overline{z})<h(z(t))$ for every $0<t<\pi/2$.
Hence $\overline{z}$ is a nonoptimal critical point of $h$. 
For $z=(z_1,z_2)^\top\in\mathcal{S}$ sufficiently close to
$\overline{z}$, write $z=z(t(z))$ with $t(z):=\arcsin(z_2)$.
Then $t(z)\to0$ as $z\to\overline{z}$ within $\mathcal{S}$,
and $t(z)>0$ whenever $h(z)>h(\overline{z})$. The above formulas for 
$h(z(t))-h(\overline{z})$ and ${\rm dist}(0,\partial h(z(t)))$ yield
\begin{equation}\label{temp-limit}
\lim_{t\downarrow0}\frac{h(z(t))-h(\overline{z})}{t^3}=\frac12\ \ {\rm and}\ \ \lim_{t\downarrow0}
\frac{{\rm dist}(0,\partial h(z(t)))}{t^2}=\frac32.
\end{equation}
For each $\theta\in[0,1)$ and $0<t<\pi/2$, write
\begin{equation}\label{Rtheta}
 R_{\theta}(t):=\frac{{\rm dist}(0,\partial h(z(t)))}{\bigl(h(z(t))-h(\overline{z})\bigr)^\theta}
=t^{2-3\theta}\frac{{\rm dist}(0,\partial h(z(t)))/t^2}
     {\bigl((h(z(t))-h(\overline{z}))/t^3\bigr)^\theta}.
\end{equation}
When $\theta=2/3$, we have $\lim_{t\downarrow0}R_{\theta}(t)=3/2^{1/3}$. 
By the local parametrization, there exists a constant $c>0$
such that ${\rm dist}(0,\partial h(z))\ge c\bigl(h(z)-h(\overline{z})\bigr)^{2/3}$ for every $z\in\mathcal{S}$ close enough to $\overline{z}$ with $h(z)>h(\overline{z})$, so $2/3$
is an admissible KL exponent of $h$ at $\overline{z}$. 

We next argue by contradiction that $2/3$ is the
smallest KL exponent of $h$ at $\overline{z}$. Fix any $\theta\in[0,2/3)$. If $\theta$ is also an admissible KL exponent of $h$ at $\overline{z}$, there would exist $c_\theta>0,\varepsilon_{\theta}>0$ and $\eta_{\theta}>0$ such that, for all $\|z-\overline{z}\|<\varepsilon_{\theta}$ and
$0<h(z)-h(\overline{z})<\eta_{\theta}$, 
\[
{\rm dist}(0,\partial h(z))
\ge c_\theta\bigl(h(z)-h(\overline{z})\bigr)^\theta.
\]
For all sufficiently small $t>0$, we have $\|z(t)-\overline{z}\|<\varepsilon_{\theta}$ and $0<h(z(t))-h(\overline{z})<\eta_{\theta}$, so the above inequality implies that $R_{\theta}(t)\ge c_{\theta}>0$. On the other hand, it follows from \eqref{temp-limit}-\eqref{Rtheta} that $\lim_{t\downarrow 0}R_{\theta}(t)=0$. Thus, we obtain a contradiction. 
\end{example}

The pointwise bound in Theorem~\ref{theorem:PLQ-34KL} can
be exact at an interface point where both adjacent quadratic
pieces have nonzero tangential gradients on the sphere.
\begin{example}\label{exa:PLQ-sharp}
Let $f$ be given by \eqref{eq:sphere-function} with
$\mathbb{E}=\mathbb{R}^3$ and $\phi(z)=(z_1-1)^2+|z_3|$ for
$z\in\mathbb{R}^3$. Set $\overline{z}:=e_1=(1,0,0)^\top$,
where $e_i$ denotes the $i$th standard basis vector of
$\mathbb{R}^3$. The function $\phi$ is PLQ, with quadratic pieces
$q_\pm(z):=(z_1-1)^2\pm z_3$ on the half-spaces $\pm z_3\ge0$.
The hyperplane $z_3=0$ induces two open half-space cells and the
interface cell $C_0:=\{z\in\mathbb R^3\mid z_3=0\}$. Both open
half-space cells have affine hull $\mathbb R^3$, so their
corresponding sphere sections are $\mathcal S$. Since
$T_{\overline z}\mathcal S=e_1^\perp$, we have
\[
\mathcal P_{T_{\overline z}\mathcal S}\nabla q_\pm(\overline z)
=\mathcal P_{e_1^\perp}(\pm e_3)=\pm e_3.
\]
Thus, by Proposition~\ref{prop:quadratic-sphere-KL}, the smallest
KL exponent of each $q_\pm+\delta_{\mathcal S}$ at $\overline z$
is $0$.

For $t\in\mathbb{R}$, let $z(t):=(\cos t,\sin t,0)^\top$. This curve
lies in $C_0\cap\mathcal{S}$ and satisfies
\[
f(z(t))-f(\overline{z})=(1-\cos t)^2\sim\tfrac14t^4
\quad\text{as }t\to0.
\]
Since $\phi$ is convex and locally Lipschitz and
$\delta_{\mathcal S}$ is subdifferentially regular on $\mathcal S$,
the subdifferential sum rule gives
$\partial f(z(t))=2(\cos t-1)e_1+[-1,1]e_3+\mathbb R z(t)$ for $t\in\mathbb{R}$. Minimizing over the normal line $\mathbb R z(t)$ amounts to projecting onto $T_{z(t)}\mathcal{S}$. Since $e_3$ is tangent and orthogonal to $e_1$, the
distance is minimized by taking the coefficient of $e_3$ to
be zero. Then ${\rm dist}(0,\partial f(z(t)))=2(1-\cos t)|\sin t|\sim |t|^3$ as $t\to0$. Consequently, 
\[
\frac{{\rm dist}(0,\partial f(z(t)))}{(f(z(t))-f(\overline z))^\theta}
\sim 4^\theta |t|^{3-4\theta}\to 0\ \ \text{as }0\ne t\to0
\]
for each $\theta\in[0,3/4)$. Hence every admissible KL exponent of $f$ at $\overline z$ is at least $3/4$. Since Theorem~\ref{theorem:PLQ-34KL} gives an admissible exponent $\theta_{\overline z}\le3/4$, the smallest KL exponent of $f$ at $\overline z$ equals
$\theta_{\overline z}=3/4$. The two open half-space cells
contribute exponent $0$, so the maximum defining
$\theta_{\overline z}$ is attained by the interface cell $C_0$.
\end{example}

Theorem~\ref{theorem:PLQ-34KL} applies to the sphere-constrained
$\ell_1$-regularized quadratic problem
\begin{equation}\label{quadratic-sphere}
\min_{z\in\mathcal S}\,\Phi(z):=\tfrac12 z^\top Qz+a^\top z+\lambda\|z\|_1,
\end{equation}
where $Q\in\mathbb R^{n\times n}$ is symmetric, $a\in\mathbb R^n$ and $\lambda\ge0$. Test problem E1 in \cite[Section~5.2]{Shen2022} has this form when the sphere radius is set to one. Specific instances include the augmented Lagrangian subproblem for finding a sparse vector in a null space and the minimization form of the $\ell_1$-regularized sparse PCA \cite[Eqs.~(5) and (7)]{Shen2022}. For each $s\in\{-1,1\}^n$, define the closed orthant $C_s:=\{z\in\mathbb R^n\mid s_jz_j\ge0\ \text{for all }j\in[n]\}$. For every $z\in C_s$, we have $|z_j|=s_jz_j$ for all $j\in[n]$, so $\|z\|_1=s^\top z$ and $\Phi(z)=\tfrac12 z^\top Qz+(a+\lambda s)^\top z$. Since the $C_s$ form a finite cover of $\mathbb{R}^n$ by closed polyhedral sets, $\Phi$ is PLQ. Hence $\Phi+\delta_{\mathcal S}$ has the KL property
with exponent $3/4$ at every point of $\mathcal S$. 

A second application is the Huber formulation of the oblique
Procrustes subproblem for a fixed target column
\cite[Eqs.~(3)--(5)]{Trendafilov2004}, which takes the form
\[
\min_{z\in\mathcal S}\,\Psi(z):=\sum_{i=1}^m\psi_\kappa(b_i^\top z-c_i),
\]
where $b_i\in\mathbb R^n$ and $c_i\in\mathbb R$, and $\psi_\kappa$ is the Huber loss with threshold $\kappa>0$ given by
\[
\psi_\kappa(t):=\left\{\begin{array}{cl}
-\kappa t-\tfrac12\kappa^2 & {\rm if}\ t\le-\kappa,\\
\tfrac12 t^2 & {\rm if}\ -\kappa\le t\le\kappa,\\
\kappa t-\tfrac12\kappa^2 &{\rm if}\ t\ge\kappa.
\end{array}\right.
\]
For each $i$, the conditions $b_i^\top z-c_i\le-\kappa$,
$-\kappa\le b_i^\top z-c_i\le\kappa$, and $b_i^\top z-c_i\ge\kappa$ define
three closed polyhedral sets that cover $\mathbb R^n$.
The nonempty intersections obtained by choosing one set for
each $i$ form a finite cover of $\mathbb R^n$ by closed polyhedral
sets. On each such intersection, every summand agrees with a polynomial
of degree at most two. Thus $\Psi$ is PLQ. Since $\psi_\kappa$ is
continuous, $\Psi$ is also continuous. Consequently,
$\Psi+\delta_{\mathcal{S}}$ satisfies the KL property with
exponent $3/4$ at every point of $\mathcal{S}$.
\subsection{KL exponents for PA functions}\label{sec3.1}

A PA function admits a PLQ representation with zero quadratic
terms. For this representation, we apply
Proposition~\ref{prop:quadratic-sphere-KL} to the affine pieces
restricted to the corresponding sections of $\mathcal{S}$ to
determine the possible values of the pointwise KL exponent
$\theta_{\overline{z}}$ in Theorem~\ref{theorem:PLQ-34KL}(i). 
\begin{theorem}\label{theorem:KL-PA}
Let $\phi$ be a proper PA function with
${\rm dom}\,\phi\cap\mathcal{S}\ne\varnothing$, and let $f$ be
given by \eqref{eq:sphere-function}. For every PLQ representation of $\phi$ with zero quadratic terms and every $\overline{z}\in{\rm dom}\,f$, the pointwise KL exponent $\theta_{\overline{z}}$ given by Theorem~\ref{theorem:PLQ-34KL}(i) lies in $\{0,1/2\}$. Consequently, $f$ satisfies the KL property with exponent $1/2$ at every point of ${\rm dom}\,f$.
\end{theorem}
\begin{proof}
Fix any PLQ representation of $\phi$ with zero quadratic
terms and use the notation of Theorem~\ref{theorem:PLQ-34KL}(i).
Fix $\overline{z}\in{\rm dom}\,f$ and
$j\in\mathcal{I}(\overline{z})$. Then
$\overline{z}\in{\rm cl}\,C_j\cap\mathcal{S}\subseteq\mathcal{M}_j$.
Set
$h_j:=q_j+\delta_{\mathcal{M}_j}$. If $\overline{z}$ is
isolated in $\mathcal{M}_j$ or $0\notin\partial h_j(\overline{z})$,
then Proposition~\ref{prop:quadratic-sphere-KL}(i) gives
$\theta_j(\overline{z})=0$. Otherwise, use the notation of
Proposition~\ref{prop:quadratic-sphere-KL}(ii) for $h_j$ at
$\overline{z}$. Since $\overline{z}$ is nonisolated in
$\mathcal{M}_j$, the tangent space
$\mathbb{V}:=T_{\overline{z}}\mathcal{M}_j$ has positive dimension.
Since $q_j$ is affine, $\mathcal{Q}=0$ and hence
$\mathcal{B}=-\mu\mathcal{I}$. Using
$\overline{u}\perp\mathbb{V}$ and $\|\overline{u}\|=1$, we obtain
\[
 \mathcal{H}=-\mu\mathcal{I}_{\mathbb{V}},\,b=0\ \ {\rm and}\ \ \gamma=-\mu,
\]
where $\mathcal{I}_{\mathbb{V}}$ is the identity mapping on $\mathbb{V}$. Since $b=0$, the vector $c_{\mathbb{L}}$ in
Lemma~\ref{lem:quartic-gradient}(ii) is zero, so
$\widehat{\gamma}=-\mu$. If $\ker\mathcal{H}\ne\{0\}$, then
$\mu=0$ and $\widehat{\gamma}=0$. Therefore,
Lemma~\ref{lem:quartic-gradient}(ii) gives exponent $1/2$ if
$\mathcal H$ has a positive eigenvalue, which occurs exactly
when $\mu<0$, and exponent $0$ otherwise. By
Proposition~\ref{prop:quadratic-sphere-KL}(ii), this is the
smallest KL exponent $\theta_j(\overline{z})$ of $h_j$ at
$\overline{z}$. Thus $\theta_j(\overline{z})\in\{0,1/2\}$ for every
$j\in\mathcal{I}(\overline{z})$, and hence
$\theta_{\overline{z}}=\max_{j\in\mathcal{I}(\overline{z})}\theta_j(\overline{z})\in\{0,1/2\}$.
By Theorem~\ref{theorem:PLQ-34KL}(i), there exist a neighborhood
$\mathcal{U}$ of $\overline{z}$ and constants $c>0$ and
$\eta\in(0,1]$ such that, for every $z\in\mathcal{U}$ with
$0<f(z)-f(\overline{z})<\eta$,
\[
\operatorname{dist}(0,\partial f(z))
\ge c\bigl(f(z)-f(\overline{z})\bigr)^{\theta_{\overline{z}}}
\ge c\bigl(f(z)-f(\overline{z})\bigr)^{1/2}.
\]
The last inequality uses $\theta_{\overline{z}}\le1/2$ and
$f(z)-f(\overline{z})<1$. Since $\overline{z}$ is arbitrary on ${\rm dom}\,f$, $f$ satisfies the KL property with exponent $1/2$ throughout ${\rm dom}\,f$.
\end{proof}
\begin{remark}\label{remark-PA1/2}
The KL exponent $1/2$ in Theorem~\ref{theorem:KL-PA} is sharp
over the PA class. Let $\mathcal{S}$ be the unit circle in
$\mathbb{R}^2$ and take $\phi(z):=1-z_1$ in
\eqref{eq:sphere-function}. The point $\overline{z}:=(1,0)^\top$
minimizes $f$ and is nonisolated in $\mathcal{S}$. In the
notation of Proposition~\ref{prop:quadratic-sphere-KL}(ii),
we have $p=0$, $r=1$, and $\overline{u}=\overline{z}$, so
$\mu=\langle\nabla\phi(\overline{z}),\overline{u}\rangle/r=-1$.
The affine-case calculation in the proof of
Theorem~\ref{theorem:KL-PA} therefore gives $1/2$ as the smallest
KL exponent of $f$ at $\overline{z}$.
\end{remark}

The KL exponents $3/4$ and $1/2$ established in 
Theorems~\ref{theorem:PLQ-34KL}(ii) and~\ref{theorem:KL-PA},
respectively, also apply when the unit sphere is intersected
with a closed polyhedral set.
\begin{corollary}\label{main-corollary}
Let $\phi:\mathbb{E}\to\overline{\mathbb{R}}$ be a proper
PLQ function and let $P\subseteq\mathbb{E}$ be a closed
polyhedral set such that ${\rm dom}\,\phi\cap P\cap\mathcal{S}\ne\varnothing$. Then $h:=\phi+\delta_{P\cap\mathcal{S}}$ satisfies the KL property
throughout ${\rm dom}\,h$ with exponent $3/4$, and with exponent
$1/2$ if $\phi$ is PA. 
\end{corollary}
\begin{proof}
Set $\widetilde{\phi}:=\phi+\delta_P$. Intersecting the closed
polyhedra in a representation of $\phi$ with $P$ shows that
$\widetilde{\phi}$ is proper PLQ and is PA whenever $\phi$ is PA.
Since ${\rm dom}\,\widetilde{\phi}\cap\mathcal{S}
\ne\varnothing$ and $h=\widetilde{\phi}+\delta_{\mathcal{S}}$,
the claims follow from Theorems~\ref{theorem:PLQ-34KL}
and~\ref{theorem:KL-PA}, respectively.
\end{proof}

As a consequence of Theorem \ref{theorem:KL-PA}, we obtain the following corollaries. 
\begin{corollary}\label{corollary1}
Let $h(z):=\max_{1\le i\le m}\langle a_i,z\rangle+\delta_{\mathcal S}(z)$ for $z\in\mathbb{E}$, where $a_i\in\mathbb{E}$ are given. Then $h$ satisfies the KL property of exponent $1/2$ at every point of $\mathcal{S}$.
\end{corollary} 
\begin{corollary}\label{corollary2-L1}
Let $A\in\mathbb R^{m\times n}$ be a given matrix, and define
$h_{+}(z):=\|Az\|_1+\delta_{\mathcal{S}}(z)$ and $h_{-}(z):=-\|Az\|_1+\delta_{\mathcal{S}}(z)$ for $z\in\mathbb{R}^n$.
Then both $h_{+}$ and $h_{-}$ satisfy the KL property of exponent $1/2$ at every point of $\mathcal{S}$. 
\end{corollary}

The pointwise maximum of finitely many linear functions and
its negative are both PA. In particular, the identity
$\|Az\|_1=\max_{s\in\{-1,1\}^m}\langle A^\top s,z\rangle$
shows that $z\mapsto\|Az\|_1$ and $z\mapsto-\|Az\|_1$ are PA.
Both corollaries follow from Theorem~\ref{theorem:KL-PA}. 
The maximum in Corollary~\ref{corollary1} is the support function $\sigma_C$ of the polytope $C\!={\rm conv}\{a_1,\ldots,a_m\}$, so the corollary applies to minimizing $\sigma_C$ over the unit sphere. Balashov and Tremba~\cite{Balashov2024-1} studied this problem for compact convex sets $C$ and gave sufficient conditions
for linear convergence of a gradient projection method. See also \cite{Balashov2024-2} for related problems. For $h_+$, Corollary~\ref{corollary2-L1} establishes the KL
property with exponent $1/2$ for arbitrary $A$, removing the
positive optimal value assumption in \cite[Theorem~4.3]{YuZengPong2025}. That assumption is equivalent
to $A$ having full column rank. The following proposition gives the
KL property with exponent $0$ for $h_+$ throughout $\mathcal S$
as a special case of a result for a broader class of
nonnegative piecewise-linear functions.
\begin{proposition}\label{prop:KL-PL}
Let $h:\mathbb{E}\to\overline{\mathbb{R}}$ be a proper
nonnegative function with ${\rm dom}\,h=\bigcup_{i=1}^m P_i$,
where each $P_i$ is a closed polyhedral cone. Suppose that for
each $i\in[m]$ there exists $a_i\in\mathbb{E}$ such that
$h(z)=\langle a_i,z\rangle$ for all $z\in P_i$.
Assume that ${\rm dom}\,h\cap\mathcal{S}\ne\varnothing$ and set
$f:=h+\delta_{\mathcal{S}}$. Then, for every
$\varpi\in f({\rm dom}\,f)$, there exists a constant $\eta_\varpi>0$
such that ${\rm dist}(0,\partial f(z))\ge\eta_{\varpi}$ for every $z\in{\rm dom}\,f$ with $0<f(z)-\varpi<\eta_{\varpi}$. Consequently, $f$ satisfies the KL property with exponent $0$
at every point of ${\rm dom}\,f$.
\end{proposition}
\begin{proof}
The function $h$ is proper PLQ with zero quadratic and constant
terms and is therefore lsc. Since $\mathcal{S}$ is closed,
$f$ is lsc. Moreover, 
${\rm dom}\,f={\rm dom}\,h\cap\mathcal{S}\ne\varnothing$,
so $f$ is proper. Lemma~\ref{lem:domain-Lsubdiff} gives
${\rm dom}\,\partial f={\rm dom}\,f$. Represent each cone $P_i$ by homogeneous linear inequalities and equalities. Use the notation of Proposition \ref{prop:tangent-proj-PLQ} with $\phi=h$ and the cells induced by the resulting hyperplanes. The cells $C_j$ for $j\in\mathcal{J}$  partition ${\rm dom}\,h$. For each $j\in\mathcal{J}$, choose $i_j\in[m]$ as in Corollary~\ref{cor:PLQ-consistency}(ii), so that $C_j\subseteq P_{i_j}$, and set $q_j(z):=\langle a_{i_j},z\rangle$ for $z\in\mathbb E$.
Since $P_{i_j}$ is closed, $q_j=h$ on ${\rm cl}\,C_j$. For each $j\in\mathcal{J}$, choose $z_j\in C_j$. The homogeneous
linear forms defining $C_j$ have the same signs at $\alpha z_j$
as at $z_j$ for every $\alpha>0$, so $\alpha z_j\in C_j$.
Letting $\alpha\downarrow0$ gives
$0\in{\rm cl}\,C_j\subseteq\mathbb{A}_j$. Thus, $\mathbb{A}_j$
is a linear subspace, and consequently, $\mathbb{A}_j=\mathbb{L}_j$
and $\mathcal{M}_j=\mathbb{L}_j\cap\mathcal{S}$. 
For each $j\in\mathcal{J}$, set $r_j:=\|\mathcal{P}_{\mathbb{L}_j}a_{i_j}\|$.

Fix any $j\in\mathcal{J}$. For any $z\in C_j\cap\mathcal{S}$, since $f(z)=q_j(z)$ and $z\in\mathbb{L}_j$, we have $f(z)=\langle\mathcal{P}_{\mathbb{L}_j}a_{i_j},z\rangle$ and $0\le f(z)\le r_j$. Since $T_z\mathcal{M}_j=\mathbb{L}_j\cap z^\perp$, Proposition~\ref{prop:tangent-proj-PLQ} gives
\[
\mathcal{P}_{T_z\mathcal{M}_j}v=\mathcal{P}_{\mathbb{L}_j}a_{i_j}-f(z)z\quad{\rm for\ every}\ \ v\in\partial f(z).
\]
The vector on the right hand side has squared norm $r_j^2-f(z)^2$. Since $\|v\|\ge\|\mathcal P_{T_z\mathcal M_j}v\|$, taking the infimum over $v\in\partial f(z)$ and using $0\le f(z)\le r_j$ yields
\begin{equation}\label{temp-KLdist}
{\rm dist}(0,\partial f(z))\ge \sqrt{r_j^2-f(z)^2}\ge r_j-f(z)
\quad\text{for every }z\in C_j\cap\mathcal{S}.
\end{equation}
Now fix any $\varpi\in f({\rm dom}\,f)$ and let
$\mathcal{J}_{\varpi}:=\{j\in\mathcal{J}\mid r_j>\varpi\}$.
If $\mathcal{J}_{\varpi}=\varnothing$, then $r_j\le\varpi$ for
every $j\in\mathcal{J}$. Since $f(z)\le r_j$ on
$C_j\cap\mathcal{S}$, we have $f(z)\le\varpi$ for every
$z\in{\rm dom}\,f$. In this case, set $\eta_{\varpi}:=1$, so there are no  $z\in{\rm dom}\,f$ satisfying $0<f(z)-\varpi<\eta_{\varpi}$. 
Otherwise, set $\eta_\varpi:=\frac12\min_{j\in\mathcal{J}_\varpi}(r_j-\varpi)>0$.
Take any $z\in{\rm dom}\,f$ with $0<f(z)-\varpi<\eta_\varpi$,
and let $j\in\mathcal{J}$ be the index such that $z\in C_j$. Since $r_j\ge f(z)>\varpi$, we have $j\in\mathcal{J}_\varpi$. The definition of $\eta_\varpi$ then gives $r_j-\varpi\ge 2\eta_\varpi$. It follows from
\eqref{temp-KLdist} that 
\[
{\rm dist}(0,\partial f(z))\ge r_j-f(z)=(r_j-\varpi)-(f(z)-\varpi)
\ge 2\eta_\varpi-(f(z)-\varpi)>\eta_\varpi.
\]
This proves the asserted estimate. For any
$\overline{z}\in{\rm dom}\,f$, applying the estimate with
$\varpi=f(\overline{z})$ gives the KL property with exponent
$0$ at $\overline{z}$.
\end{proof}

The support function $\sigma_C$ of a convex polytope $C$
containing the origin is nonnegative and linear on each of
finitely many closed polyhedral cones that cover $\mathbb{E}$.
We apply Proposition~\ref{prop:KL-PL} with $h=\sigma_C$ and
use the resulting subdifferential distance estimate to prove
finite convergence of the full-step ManPPA iteration for
minimizing $\sigma_C$ over $\mathcal{S}$. This iteration is
adapted from \cite[Algorithm~1]{ChenSo2021}.

\begin{corollary}\label{cor:manppa-finite}
Let $C\subseteq\mathbb{E}$ be a convex polytope containing the
origin, and set 
\[
 h:=\sigma_C\ \ {\rm and}\ \ f:=h+\delta_{\mathcal{S}},
\]
where $\sigma_C(z):=\max_{u\in C}\langle u,z\rangle$.
Fix any $t>0$ and $z^0\in\mathcal{S}$. For each $k\ge0$, let
\begin{equation}\label{manPPA-iteration}
 d^k:=\mathop{\arg\min}_{d\in(z^k)^\perp}
\Big\{h(z^k+d)+\frac{1}{2t}\|d\|^2\Big\}\ \ {\rm and}\ \ 
z^{k+1}:=\frac{z^k+d^k}{\|z^k+d^k\|}.
\end{equation}
Then there exists an integer $K\ge0$ such that $z^k=z^K$ for
all $k\ge K$ and $0\in\partial f(z^K)$.
\end{corollary}
\begin{proof}
Let $a_1,\ldots,a_m\in\mathbb{E}$ be the vertices of $C$. Then
$C={\rm conv}\{a_1,\ldots,a_m\}$, and hence
$h(z)=\max_{i\in[m]}\langle a_i,z\rangle$ for every $z\in\mathbb{E}$. For each $i\in[m]$, define
\[
 P_i:=\big\{z\in\mathbb{E}\mid
\langle a_i-a_j,z\rangle\ge0\ \ \text{for every }j\in[m]\big\}.
\]
The closed polyhedral cones $P_i$ cover $\mathbb{E}$, and
$h(z)=\langle a_i,z\rangle$ for all $z\in P_i$. Since $0\in C$, we have $h\ge 0$. Thus $h$ satisfies the assumptions of Proposition~\ref{prop:KL-PL}. By \cite[Corollary~23.5.3]{Roc70}, the subdifferential of $h$ takes the following form 
\[
\partial h(y)=\big\{u\in C\mid \langle u,y\rangle=h(y)\big\}
\quad\text{for every }y\in\mathbb{E}.
\]
Then every subgradient of $h$ has norm at most $L:=\max_{u\in C}\|u\|$, and $\partial h(\alpha y)=\partial h(y)$ for every $\alpha>0$. Since $h$ is convex, $\widehat{\partial}h(z)=\partial h(z)$ for every $z\in\mathbb{E}$. The normal-space
formula for the sphere gives 
$\widehat{\partial}\delta_{\mathcal{S}}(z)=\mathbb{R}z$
for every $z\in\mathcal{S}$. Hence the addition rule in
\cite[Corollary~10.9]{RW98} yields $\partial h(z)+\mathbb{R}z\subseteq\widehat{\partial}f(z)
\subseteq\partial f(z)$ for every $z\in\mathcal{S}$.

Each subproblem in \eqref{manPPA-iteration} has a unique minimizer because its
objective is continuous, coercive, and strongly convex on
$(z^k)^\perp$. Set $y^k:=z^k+d^k$. Since $d^k\perp z^k$,
we have $\|y^k\|=\sqrt{1+\|d^k\|^2}\ge 1$. Thus the normalization is
well defined and $z^{k+1}\in\mathcal{S}$. By nonnegativity
and positive homogeneity of $h$, we have 
$h(z^{k+1})=h(y^k)/\|y^k\|\le h(y^k)$. Note that $d=0$ is feasible for every subproblem. The optimality of $d^k$ yields that
\begin{equation}\label{eq:manppa-sphere-descent} 
f(z^{k+1})+\frac{1}{2t}\|d^k\|^2\le f(z^k)\quad{\rm for\ every}\ k.
\end{equation}
By \eqref{eq:manppa-sphere-descent} and the nonnegativity of
$f$, we have $f(z^k)\downarrow\varpi^*$ for some $\varpi^*\ge0$.
Summing \eqref{eq:manppa-sphere-descent} over $k$ leads to 
$\sum_{k=0}^{\infty}\|d^k\|^2\le 2t\bigl(f(z^0)-\varpi^*\bigr)<\infty$. Hence $d^k\to0$. By the optimality condition of subproblems, there exist
$v^k\in\partial h(y^k)$ and $\lambda_k\in\mathbb{R}$ such that
\[
 v^k+t^{-1}d^k+\lambda_k z^k=0. 
\]
Taking the inner product with $z^k$ yields $|\lambda_k|=|\langle v^k,z^k\rangle|\le L$.
Since $\partial h(\alpha y)=\partial h(y)$ for every $\alpha>0$, and $z^{k+1}=y^k/\|y^k\|$, we have $v^k\in\partial h(z^{k+1})$. As $y^k$ is a positive multiple of $z^{k+1}$, it follows from the above equality that
\begin{equation}\label{zk-inclusion}
(\lambda_k-t^{-1})d^k
=v^k+\lambda_k y^k\in\partial f(z^{k+1}).
\end{equation}
Consequently, ${\rm dist}(0,\partial f(z^{k+1}))\le (t^{-1}+L)\|d^k\|\to 0$. 
Since $h$ is continuous and $\mathcal{S}$ is compact, $f({\rm dom}\,f)=h(\mathcal{S})$ is compact. Thus
$\varpi^*\in f({\rm dom}\,f)$. Proposition~\ref{prop:KL-PL} 
gives a constant $\eta_{\varpi^*}>0$ such that, for every $z\in{\rm dom}\,f$ with $0<f(z)-\varpi^*<\eta_{\varpi^*}$,  
\begin{equation}\label{temp-ineqKL}
{\rm dist}(0,\partial f(z))\ge\eta_{\varpi^*}.
\end{equation} 
Since $f(z^k)\downarrow\varpi^*$
and ${\rm dist}(0,\partial f(z^k))\to0$, we can choose
$K\ge1$ such that 
\[
 0\le f(z^K)-\varpi^*<\eta_{\varpi^*}\ \ {\rm and}\ \ {\rm dist}(0,\partial f(z^K))<\eta_{\varpi^*}. 
\]
If $f(z^K)>\varpi^*$, then $0<f(z^K)-\varpi^*<\eta_{\varpi^*}$, so
\eqref{temp-ineqKL} gives ${\rm dist}(0,\partial f(z^K))\ge\eta_{\varpi^*}$,
contradicting the choice of $K$. Hence $f(z^K)=\varpi^*$. By monotonicity, $f(z^k)=\varpi^*$ for all $k\ge K$. It follows from \eqref{eq:manppa-sphere-descent} that $d^k=0$ for all $k\ge K$, so \eqref{manPPA-iteration}
gives $z^k=z^K$ for all $k\ge K$. Taking $k=K$ in \eqref{zk-inclusion} yields
$0\in\partial f(z^{K+1})=\partial f(z^K)$. 
\end{proof}
\begin{remark}\label{remark-manppa}
{\bf(a)} The iteration in Corollary~\ref{cor:manppa-finite} can also
be implemented using the backtracking line search in
\cite[Algorithm~1]{ChenSo2021} to select the stepsize
$\alpha_k$ along the tangent direction $d^k$. Fix any $\beta\in(0,1)$. Given a current iterate $z^k$ and the minimizer $d^k$ of the tangent subproblem, choose $\alpha_k$ as the first trial $\alpha$ in $1,\beta,\beta^2,\ldots$ satisfying
\[
h\left(\frac{z^k+\alpha d^k}{\|z^k+\alpha d^k\|}\right)
\le h(z^k)-\frac{\alpha}{2t}\|d^k\|^2,
\]
and set $z^{k+1}:=(z^k+\alpha_k d^k)/\|z^k+\alpha_k d^k\|$. Since $f=h$ on $\mathcal{S}$, the descent estimate
\eqref{eq:manppa-sphere-descent} shows that $\alpha=1$ satisfies this condition at every iteration. Hence $\alpha_k=1$ for all $k\ge0$. The backtracking iteration thus generates the same iterates as the full-step iteration and stabilizes at a critical point of $f$
in finitely many steps.

\noindent
\textbf{(b)} For $A\in\mathbb{R}^{m\times n}$, the polytope
$C=A^\top[-1,1]^m$ contains the origin and satisfies
\[
\sigma_C(z)=\max_{s\in[-1,1]^m}s^\top Az=\|Az\|_1.
\]
For this choice of $C$, the backtracking iteration in part~(a)
coincides with Algorithm~1 of \cite{ChenSo2021} when the tangent
subproblems are solved exactly. By Corollary~\ref{cor:manppa-finite}
and part~(a), its iterates stabilize at a critical point in
finitely many steps for every $A$, every $z^0\in\mathcal{S}$,
and every fixed proximal parameter $t>0$. This finite convergence guarantee removes the assumption in \cite[Theorem~4.8]{YuZengPong2025} that the limit
lies in a set of weak sharp minima. 

For the same sphere-constrained $\ell_1$-norm minimization problem,
Tsakiris and Vidal~\cite[Proposition~16]{TsakirisVidal2018} gave
a proof that the recursive linear-programming method of
Sp\"ath and Watson~\cite{SpathWatson1987} reaches a critical
point in finitely many steps. Their result assumes that $A$
has full column rank and that the solution selected for each
linear program is orthogonal to $n-1$ linearly independent
rows of $A$.
\end{remark}

\subsection{Homogeneous quadratic pieces on polyhedral cones}\label{sec3.3}

Theorem~\ref{theorem:PLQ-34KL} establishes the KL property with
exponent $3/4$ for $f$ throughout ${\rm dom}\,f$. The following
proposition shows that the exponent can be reduced to $1/2$
when ${\rm dom}\,\phi$ is a finite union of closed convex
polyhedral cones and $\phi$ is homogeneous quadratic on each
cone. Its proof uses spectral estimates for the quadratic forms
restricted to the linear spans of the cone faces.
\begin{proposition}\label{prop:ph12}
Let $\phi:\mathbb{E}\to\overline{\mathbb{R}}$ be a proper
lsc function with ${\rm dom}\,\phi=\bigcup_{i=1}^m K_i$, where each $K_i\subseteq\mathbb{E}$ is a closed convex polyhedral cone. Suppose that for each $i\in[m]$ there exists a self-adjoint linear operator $\mathcal{Q}_i:\mathbb{E}\to\mathbb{E}$ such that $\phi(z)=\frac12\langle z,\mathcal{Q}_i z\rangle$ for all $z\in K_i$. Assume that ${\rm dom}\,\phi\cap\mathcal{S}\ne\varnothing$ 
and let $f$ be given by \eqref{eq:sphere-function}. Then, for every $\varpi\in f({\rm dom}\,f)$, there exists $\gamma_{\varpi}>0$ such that
${\rm dist}(0,\partial f(z))\ge \sqrt{\gamma_{\varpi}(f(z)-\varpi)}$ for every $z\in{\rm dom}\,f$ satisfying
$0<f(z)-\varpi<\gamma_{\varpi}/4$. Consequently, $f$ satisfies the KL property with exponent $1/2$ at every point of ${\rm dom}\,f$. 
\end{proposition}
\begin{proof}
The function $f$ is proper and lsc, with ${\rm dom}\,f=\bigcup_{i=1}^m(K_i\cap\mathcal{S})$. Since $\phi$ satisfies the assumptions of Corollary~\ref{cor:PLQ-consistency} with $P_i=K_i$ and zero linear and constant terms, Lemma~\ref{lem:domain-Lsubdiff} gives ${\rm dom}\,\partial f={\rm dom}\,f$. Fix any $\varpi\in f({\rm dom}\,f)$. For each $i\in[m]$ and each face $F$ of $K_i$ (including $K_i$ itself) with $F\cap\mathcal{S}\ne\varnothing$, define
\[
 \mathbb{L}_F:={\rm span}\,F\ \ {\rm and}\ \  \mathcal{A}_{i,F}:=\mathcal{P}_{\mathbb{L}_F}\mathcal{Q}_i|_{\mathbb{L}_F}.
\]
For $u,v\in\mathbb{L}_F$, the self-adjointness of 
$\mathcal{Q}_i$ gives $\langle u,\mathcal{A}_{i,F}v\rangle=\langle u,\mathcal{Q}_i v\rangle=\langle\mathcal{Q}_iu,v\rangle=\langle\mathcal{A}_{i,F}u,v\rangle$. Thus $\mathcal{A}_{i,F}$ is self-adjoint on $\mathbb{L}_F$.
Let $\sigma(\mathcal{A}_{i,F})$ denote its set of eigenvalues. Since each $K_i$ has finitely many faces and each operator $\mathcal{A}_{i,F}$ has a finite spectrum, we can choose
$\gamma_{\varpi}>0$ such that, for each $i\in[m]$ and each face $F$
of $K_i$ with $F\cap\mathcal{S}\ne\varnothing$, 
\[
 \mu-2\varpi\ge\gamma_{\varpi}\quad{\rm for\ every}\ \mu\in\sigma(\mathcal{A}_{i,F})\cap(2\varpi,\infty).
\]
Then every $\mu\in\sigma(\mathcal{A}_{i,F})$ satisfies
$(\mu-2\varpi)^2\ge\gamma_{\varpi}(\mu-2\varpi)$. By the spectral theorem, $\mathbb{L}_F$ admits an orthonormal basis consisting of eigenvectors of $\mathcal{A}_{i,F}$. Consequently, for each $u\in\mathbb{L}_F$, expanding $u$ in this orthonormal basis and applying the preceding scalar inequality to every eigenvalue yields
\begin{equation}\label{eq:ph-spectral}
\|(\mathcal{A}_{i,F}-2\varpi\,\mathcal{I}_{\mathbb{L}_F})u\|^2
\ge\gamma_{\varpi}\langle u,(\mathcal{A}_{i,F}-2\varpi\,\mathcal{I}_{\mathbb{L}_F})u\rangle.
\end{equation}

Fix any $z\in\mathbb{E}$ with $0<f(z)-\varpi<\gamma_{\varpi}/4$. Clearly, $z\in{\rm dom}\,f$. Choose $i\in[m]$ such that $z\in K_i$ and let $F$ be the minimal face of $K_i$ containing $z$. Then $z\in{\rm ri}(F)$ and $z\in F\cap\mathcal{S}$. Set $T:=\mathbb{L}_F\cap z^\perp$ and $\mathcal{A}:=\mathcal{A}_{i,F}$. For each $d\in T\setminus\{0\}$, define 
\[
 w(t):=\frac{z+td}{\|z+td\|}=\frac{z+td}{\sqrt{1+t^2\|d\|^2}}\quad{\rm for}\ t\in\mathbb{R}.
\]
Then $w(0)=z$ and $w'(0)=d$. Since $z\in{\rm ri}(F)$ and
${\rm aff}(F)=\mathbb{L}_F$, we have $z+td\in F$ for all
sufficiently small $|t|$. As $F$ is a cone, we have $w(t)\in F\cap\mathcal{S}$ for these $t$, and 
\[
 f(w(t))=\phi(w(t))=\tfrac12\langle w(t),\mathcal{Q}_iw(t)\rangle=\tfrac12\langle w(t),\mathcal{A}w(t)\rangle. 
\]
Take any $\xi\in\widehat{\partial}f(z)$. For each $d\in T\setminus\{0\}$, substituting $w(t)=z+td+o(|t|)$ into the quadratic expression for $f(w(t))$ 
and expanding gives 
\[
\frac{f(w(t))-f(z)-\langle\xi,w(t)-z\rangle}{\|w(t)-z\|}
=\frac{t}{|t|}\frac{\langle\mathcal{A}z-\xi,d\rangle}{\|d\|}
+o(1)\quad\text{as }0\ne t\to0.
\]
By the definition of Fr\'echet subgradients, both one-sided
limits are nonnegative, so
$\langle\xi,d\rangle=\langle\mathcal{A}z,d\rangle$.
This equality holds for $d=0$, so for all $d\in T$. Then $\mathcal{P}_T\xi=\mathcal{P}_T\mathcal{A}z$. 
Set $\rho:=\langle z,\mathcal{A}z\rangle=2f(z)$. Since
$\mathcal{A}z\in\mathbb{L}_F$ and $\|z\|=1$, we have
$\mathcal{P}_T\xi=\mathcal{P}_T\mathcal{A}z=(\mathcal{A}-\rho\mathcal{I}_{\mathbb{L}_F})z$. 
Since $(\mathcal{A}-\rho\mathcal{I}_{\mathbb{L}_F})z$ is orthogonal to $z$, using $\|z\|=1$ and \eqref{eq:ph-spectral} yields
\begin{align}\label{Rsubdiff-bound}
\|\xi\|^2\ge\|\mathcal{P}_T\xi\|^2&=\|(\mathcal{A}-\rho\mathcal{I}_{\mathbb{L}_F})z\|^2=\|(\mathcal{A}-2\varpi\mathcal{I}_{\mathbb{L}_F})z\|^2-(\rho-2\varpi)^2\nonumber\\
&\ge(\rho-2\varpi)\bigl(\gamma_{\varpi}-(\rho-2\varpi)\bigr)\ge\gamma_{\varpi}\bigl(f(z)-\varpi\bigr),
\end{align}
where the last inequality uses $0<\rho-2\varpi=2(f(z)-\varpi)<\gamma_{\varpi}/2$. 

Fix any $z\in\mathbb{E}$ with $0<f(z)-\varpi<\gamma_{\varpi}/4$ and any $v\in\partial f(z)$. There exist sequences $z^k\to z$ and $v^k\to v$ such that $f(z^k)\to f(z)$, and $v^k\in\widehat{\partial}f(z^k)$ for each $k\in\mathbb{N}$. Since $0<f(z)-\varpi<\gamma_{\varpi}/4$, we have $f(z^k)\in(\varpi,\varpi+\gamma_{\varpi}/4)$ for all sufficiently large $k$. Applying the above \eqref{Rsubdiff-bound} with $(z,\xi)=(z^k,v^k)$ for sufficiently large $k$ gives $\|v^k\|^2\ge\gamma_{\varpi}(f(z^k)-\varpi)$. Passing to the limit in this inequality, taking square roots, and then taking the infimum over $v\in\partial f(z)$ yields ${\rm dist}(0,\partial f(z))
\ge\sqrt{\gamma_{\varpi}(f(z)-\varpi)}$. This proves the asserted estimate for every $\varpi\in f({\rm dom}\,f)$. For the second conclusion, fix any $\overline{z}\in{\rm dom}\,f$. Applying the estimate with $\varpi=f(\overline{z})$ and taking $\eta=\gamma_{\varpi}/4$ gives the KL property with exponent $1/2$ at $\overline{z}$.
\end{proof}
\begin{remark}\label{rem:applications-ph12}
Proposition~\ref{prop:ph12} provides the KL property with exponent
$1/2$ for the extended-valued objectives of two sphere-constrained
problems, namely the rank-one sparse PCA model with a squared
$\ell_1$-norm penalty and the minimization problem defining the second
variational value of a hypergraph Laplacian. For sparse PCA, impose the rank-one constraint in the penalized semidefinite program of \cite[Eq.~(3.2)]{dAspremont2007} by
writing $X=xx^\top$ with $x\in\mathbb{R}^n$. Then, ${\rm tr}(X)=\|x\|^2$, ${\rm tr}(AX)=x^\top A x$, and 
$\rho\sum_{j,k}|X_{jk}|=\rho\|x\|_1^2$. Negating the maximization
objective gives
\[
\min_{x\in\mathcal{S}}\,\phi_{\mathrm{PCA}}(x):=\rho\|x\|_1^2-x^\top A x,
\]
where $\rho>0$ is the penalty parameter and $A$ is the covariance matrix. For each $s\in\{-1,1\}^n$, define $C_s\!:=\{x\in\mathbb{R}^n\mid s_jx_j\ge0\ \text{for all }j\in[n]\}$. These closed convex polyhedral cones cover 
${\rm dom}\,\phi_{\mathrm{PCA}}=\mathbb{R}^n$.
On each $C_s$, we have $\|x\|_1=s^\top x$ and hence
$\phi_{\mathrm{PCA}}(x)=\tfrac12\langle x,\mathcal Q_s x\rangle$, where $\mathcal{Q}_s:=2(\rho ss^\top\!-\!A)$ is symmetric.  

For the hypergraph Laplacian in \cite{Chan2018}, let $(V,E)$
be a weighted hypergraph, where $V$ is the vertex set with
$|V|=n\ge2$ and $E$ is the set of hyperedges. Each hyperedge
$e\in E$ is a nonempty subset of $V$ and has weight $w_e>0$.
Assume that the weighted degree
$d_v:=\sum_{e\in E,\,v\in e}w_e$ is positive for every $v\in V$.
Set $W:={\rm diag}(d_v)_{v\in V}$ and let $\mathbf{1}\in\mathbb{R}^n$
denote the all-ones vector. Normalizing the denominator of the
Rayleigh quotient to one gives the following problem for the
second variational value:
\[
\min_{\substack{z^\top Wz=1\\ z^\top W\mathbf{1}=0}} h(z):=\sum_{e\in E}w_e
\Big(\max_{v\in e}z_v-\min_{v\in e}z_v\Big)^2.
\]
Under the change of variables $x=W^{1/2}z$, the feasible set
becomes $\mathcal{S}\cap\mathbb{L}$, where
$\mathbb{L}:=(W^{1/2}\mathbf 1)^\perp$. Thus the problem can be written as $\min_{x\in\mathcal{S}}\,\phi_{\mathrm H}(x):=h(W^{-1/2}x)+\delta_{\mathbb{L}}(x)$. For each ordering $\pi$ of the vertices, let
$C_\pi:=\big\{x\in \mathbb{L}\mid
(W^{-1/2}x)_{\pi(1)}\le\cdots\le(W^{-1/2}x)_{\pi(n)}\big\}$. 
These closed polyhedral cones cover $\mathbb{L}$. For each hyperedge $e$, its first and last vertices in the ordering $\pi$ attain the
minimum and maximum, respectively, of $(W^{-1/2}x)_v$ over
$v\in e$ for every $x\in C_\pi$. Hence $\phi_{\mathrm H}$ is
a weighted sum of squares of linear forms on each $C_\pi$.
Since $\dim L=n-1\ge1$, we have $\mathcal{S}\cap\mathbb{L}\ne\varnothing$.

Both $\phi_{\mathrm{PCA}}$ and $\phi_{\mathrm H}$ are proper
and lsc. Proposition~\ref{prop:ph12} therefore
gives the KL property with exponent $1/2$ for
$\phi_{\mathrm{PCA}}+\delta_{\mathcal S}$ on $\mathcal S$ and for
$\phi_{\mathrm H}+\delta_{\mathcal S}$ on $\mathcal S\cap\mathbb{L}$.
\end{remark}

Example~\ref{exa:PLQ-sharp} also shows that the smallest KL
exponent on $\mathcal{S}$ can increase when a homogeneous
quadratic is added to a piecewise-affine function.
In that example, $\phi=\vartheta+q$, where
$\vartheta(z):=1-2z_1+|z_3|$ and $q(z):=z_1^2$.
By Theorem~\ref{theorem:KL-PA} and Proposition~\ref{prop:ph12},
both $\vartheta+\delta_{\mathcal{S}}$ and
$q+\delta_{\mathcal{S}}$ admit KL exponent $1/2$ at
$\overline{z}=e_1$. Their sum is $f$, whose smallest KL
exponent at $\overline{z}$ is $3/4$.

\section{Conclusion}\label{sec4}

We developed a systematic analysis of explicit KL exponents
for proper PLQ functions restricted to the unit sphere.
By combining the classification of the smallest KL exponents
of quadratics restricted to affine sections of $\mathcal{S}$
with a tangential projection identity, we obtained an admissible
exponent $\theta_{\overline{z}}$ for $f$ at each reference point
and a uniform exponent of $3/4$ that is independent of the
ambient dimension and sharp over the PLQ class. At the nonsmooth
interface in Example~\ref{exa:PLQ-sharp}, the smallest KL
exponent of $f$ is $3/4$, and each adjacent quadratic piece
restricted to the full sphere has smallest KL exponent $0$
at the same point. For proper PA $\phi$, we obtained the sharp uniform
exponent $1/2$. When ${\rm dom}\,\phi$ is a finite union of
closed polyhedral cones, $f$ admits exponent $1/2$ throughout
its domain if $\phi$ is homogeneous quadratic on each cone,
and exponent $0$ if $\phi$ is nonnegative and linear on each
cone. 

For feasible sequences converging to $\overline{z}$ and satisfying the sufficient decrease and relative error conditions with fixed positive constants, the uniform KL exponent $3/4$ yields the function-value estimate
$f(z^k)-f(\overline{z})=O(k^{-2})$ by \cite[Theorem~3.4]{Frankel2015}. When
$\theta_{\overline{z}}=2/3$, the pointwise exponent sharpens
this estimate to $O(k^{-3})$. These estimates also apply when
$\overline{z}$ is a nonoptimal critical point. For any convex
polytope $C$ containing the origin, we used exponent $0$ to
prove that full-step ManPPA for minimizing $\sigma_C$ over
$\mathcal{S}$ stabilizes at a critical point in finitely many
steps from any feasible initial point and for any fixed
positive proximal parameter.

\section*{Acknowledgments}
The authors used ChatGPT to help develop and check
parts of the proofs and to improve the exposition.
The authors assume responsibility for all content.

\bibliography{reference}\bibliographystyle{plain}

@article{dAspremont2007,
  author  = {d'Aspremont, Alexandre and El Ghaoui, Laurent and
             Jordan, Michael I. and Lanckriet, Gert R. G.},
  title   = {A direct formulation for sparse {PCA} using
             semidefinite programming},
  journal = {SIAM Review},
  year    = {2007},
  volume  = {49},
  number  = {3},
  pages   = {434--448},
  doi     = {10.1137/050645506}
}

@article{Attouch2013,
  author  = {Attouch, Hedy and Bolte, J{\'e}r{\^o}me and
             Svaiter, Benar Fux},
  title   = {Convergence of descent methods for semi-algebraic and
             tame problems: proximal algorithms, forward--backward
             splitting, and regularized {Gauss--Seidel} methods},
  journal = {Mathematical Programming},
  year    = {2013},
  volume  = {137},
  number  = {1--2},
  pages   = {91--129},
  doi     = {10.1007/s10107-011-0484-9}
}

@article{Balashov2024-1,
  author  = {Balashov, M. V. and Tremba, A. A.},
  title   = {The Gradient Projection Method for a Supporting Function on the Unit Sphere and Its Applications},
  journal = {Computational Mathematics and Mathematical Physics},
  year    = {2024},
  volume  = {64},
  number  = {4},
  pages   = {676--692},
  doi     = {10.1134/S096554252470009X}
}

@article{Balashov2024-2,
  author  = {Balashov, M. V. and Biglov, K. Z. and Tremba, A. A.},
  title   = {On Some Problems with Multivalued Mappings},
  journal = {Automation and Remote Control},
  year    = {2024},
  volume  = {85},
  number  = {5},
  pages   = {491--511},
  doi     = {10.1134/S0005117924050035}
}

@article{Bolte2007,
  author  = {Bolte, J{\'e}r{\^o}me and Daniilidis, Aris and Lewis, Adrian},
  title   = {The {L}ojasiewicz Inequality for Nonsmooth Subanalytic
             Functions with Applications to Subgradient Dynamical Systems},
  journal = {SIAM Journal on Optimization},
  volume  = {17},
  number  = {4},
  pages   = {1205--1223},
  year    = {2007},
  doi     = {10.1137/050644641}
}

@article{Bolte2007-II,
  author  = {Bolte, J{\'e}r{\^o}me and Daniilidis, Aris and
             Lewis, Adrian S. and Shiota, Masahiro},
  title   = {Clarke Subgradients of Stratifiable Functions},
  journal = {SIAM Journal on Optimization},
  volume  = {18},
  number  = {2},
  pages   = {556--572},
  year    = {2007},
  doi     = {10.1137/060670080}
}

@article{Chan2018,
  author  = {Chan, T.-H. Hubert and Louis, Anand and
             Tang, Zhihao Gavin and Zhang, Chenzi},
  title   = {Spectral properties of hypergraph {Laplacian}
             and approximation algorithms},
  journal = {Journal of the ACM},
  year    = {2018},
  volume  = {65},
  number  = {3},
  pages   = {15:1--15:48},
  doi     = {10.1145/3178123}
}

@article{LiPong2018,
  author  = {Li, Guoyin and Pong, Ting Kei},
  title   = {Calculus of the Exponent of {Kurdyka--Lojasiewicz}
             Inequality and Its Applications to Linear Convergence
             of First-Order Methods},
  journal = {Foundations of Computational Mathematics},
  volume  = {18},
  number  = {5},
  pages   = {1199--1232},
  year    = {2018},
  doi     = {10.1007/s10208-017-9366-8}
}

@article{QuSunWright2016,
  title   = {Finding a Sparse Vector in a Subspace: Linear Sparsity Using Alternating Directions},
  author  = {Qu, Qing and Sun, Ju and Wright, John},
  journal = {IEEE Transactions on Information Theory},
  volume  = {62},
  number  = {10},
  pages   = {5855--5880},
  year    = {2016},
  month   = oct,
  doi     = {10.1109/TIT.2016.2601599}
}

@article{Frankel2015,
  author  = {Frankel, Pierre and Garrigos, Guillaume and Peypouquet, Juan},
  title   = {Splitting methods with variable metric for
             {Kurdyka--{\L}ojasiewicz} functions and general convergence rates},
  journal = {Journal of Optimization Theory and Applications},
  year    = {2015},
  volume  = {165},
  number  = {3},
  pages   = {874--900},
  doi     = {10.1007/s10957-014-0642-3}
}

@article{JiangDuanLi2019,
  author  = {Jiang, Rujun and Li, Duan},
  title   = {Novel Reformulations and Efficient Algorithms for the
             Generalized Trust Region Subproblem},
  journal = {SIAM Journal on Optimization},
  year    = {2019},
  volume  = {29},
  number  = {2},
  pages   = {1603--1633},
  doi     = {10.1137/18M1174313}
}

@article{JiangLi2022,
  author  = {Jiang, Rujun and Li, Xudong},
  title   = {H{\"o}lderian error bounds and {Kurdyka-{\L}ojasiewicz}
             inequality for the trust region subproblem},
  journal = {Mathematics of Operations Research},
  year    = {2022},
  volume  = {47},
  number  = {4},
  pages   = {3025--3050},
  doi     = {10.1287/moor.2021.1243}
}

@article{Trendafilov2004,
  author  = {Trendafilov, Nickolay T. and Watson, G. A.},
  title   = {The {$\ell_1$} oblique {Procrustes} problem},
  journal = {Statistics and Computing},
  year    = {2004},
  volume  = {14},
  number  = {1},
  pages   = {39--51},
  doi     = {10.1023/B:STCO.0000009415.14785.2a}
}

@article{Kwak2008,
  title   = {Principal Component Analysis Based on {$L_1$}-Norm Maximization},
  author  = {Kwak, Nojun},
  journal = {IEEE Transactions on Pattern Analysis and Machine Intelligence},
  volume  = {30},
  number  = {9},
  pages   = {1672--1680},
  year    = {2008},
  month   = sep,
  doi     = {10.1109/TPAMI.2008.114}
}

@article{TsakirisVidal2018,
  title   = {Dual Principal Component Pursuit},
  author  = {Tsakiris, Manolis C. and Vidal, Ren{\'e}},
  journal = {Journal of Machine Learning Research},
  volume  = {19},
  number  = {18},
  pages   = {1--50},
  year    = {2018},
  url     = {https://jmlr.org/papers/v19/17-436.html}
}

@article{WuPanBi2021,
  title   = {Kurdyka--{\L}ojasiewicz Property of Zero-Norm Composite Functions},
  author  = {Wu, Yuqia and Pan, Shaohua and Bi, Shujun},
  journal = {Journal of Optimization Theory and Applications},
  volume  = {188},
  number  = {1},
  pages   = {94--112},
  year    = {2021},
  doi     = {10.1007/s10957-020-01779-7}
}

@article{YuLiPong2022,
  author  = {Yu, Peiran and Li, Guoyin and Pong, Ting Kei},
  title   = {{Kurdyka--{\L}ojasiewicz} exponent via inf-projection},
  journal = {Foundations of Computational Mathematics},
  year    = {2022},
  volume  = {22},
  pages   = {1171--1217},
  doi     = {10.1007/s10208-021-09528-6}
}

@article{Shen2022,
  author  = {Shen, Chungen and Mi, Ling and Zhang, Lei-Hong},
  title   = {An active-set proximal quasi-Newton algorithm for
             $\ell_1$-regularized minimization over a sphere constraint},
  journal = {Optimization},
  volume  = {71},
  number  = {16},
  pages   = {4623--4664},
  year    = {2022},
  doi     = {10.1080/02331934.2021.1958809}
}

@article{Montanari2016,
  author  = {Montanari, Andrea and Richard, Emile},
  title   = {Non-Negative Principal Component Analysis:
             Message Passing Algorithms and Sharp Asymptotics},
  journal = {IEEE Transactions on Information Theory},
  volume  = {62},
  number  = {3},
  pages   = {1458--1484},
  year    = {2016},
  doi     = {10.1109/TIT.2015.2457942}
}

@article{Ferreira2024,
  author  = {Ferreira, Orizon Pereira and
             Gao, Yingchao and
             N{\'e}meth, S{\'a}ndor Zolt{\'a}n and
             Rig{\'o}, Petra Ren{\'a}ta},
  title   = {Gradient Projection Method on the Sphere,
             Complementarity Problems and Copositivity},
  journal = {Journal of Global Optimization},
  volume  = {90},
  number  = {1},
  pages   = {1--25},
  year    = {2024},
  doi     = {10.1007/s10898-024-01390-4}
}

@article{YuZengPong2025,
  author  = {Peiran Yu and Liaoyuan Zeng and Ting Kei Pong},
  title   = {Convergence Analysis for a Variant of Manifold Proximal
             Point Algorithm Based on Kurdyka--{\L}ojasiewicz Property},
  journal = {Pacific Journal of Optimization},
  volume  = {21},
  pages   = {761--785},
  year    = {2025},
  doi     = {10.61208/pjo-2025-037}
}

@article{LiuSo2019,
  author  = {Huikang Liu and Anthony Man-Cho So and Weijie Wu},
  title   = {Quadratic Optimization with Orthogonality Constraint:
             Explicit {\L}ojasiewicz Exponent and Linear Convergence
             of Retraction-Based Line-Search and Stochastic
             Variance-Reduced Gradient Methods},
  journal = {Mathematical Programming},
  volume  = {178},
  number  = {1--2},
  pages   = {215--262},
  year    = {2019},
  doi     = {10.1007/s10107-018-1285-1}
}

@article{ChenSo2021,
  author  = {Chen, Shixiang and Deng, Zengde and Ma, Shiqian
             and So, Anthony Man-Cho},
  title   = {Manifold proximal point algorithms for dual principal
             component pursuit and orthogonal dictionary learning},
  journal = {IEEE Trans. Signal Process.},
  volume  = {69},
  year    = {2021},
  pages   = {4759--4773},
  doi     = {10.1109/TSP.2021.3099643}
}

@book{RW98,
	title={{Variational Analysis}},
	author={Rockafellar, R Tyrrell and Wets, Roger JB},
	year={1998},
	publisher={Springer}
}

@article{SpathWatson1987,
  author  = {Sp{\"a}th, H. and Watson, G. A.},
  title   = {On orthogonal linear {$\ell_1$} approximation},
  journal = {Numer. Math.},
  year    = {1987},
  volume  = {51},
  pages   = {531--543},
  doi     = {10.1007/BF01400354}
}

@book{Roc70,
 title={{Convex Analysis}},
 author={Rockafellar, R Tyrrell},
 year={1970},
 publisher={Princeton university press}
}

@misc{WuLeeXie2026,
  author = {Wu, Hanju and Lee, Ching-pei and Xie, Yue},
  title  = {Transferability of error bounds and the
            {Kurdyka}--{{\L}ojasiewicz} property under {$C^1$}
            partial smoothness},
  year   = {2026},
  note   = {arXiv:2605.02754v4}
}

@article{DAcuntoKurdyka2005,
  author  = {D'Acunto, Didier and Kurdyka, Krzysztof},
  title   = {Explicit bounds for the {{\L}ojasiewicz} exponent
             in the gradient inequality for polynomials},
  journal = {Ann. Polon. Math.},
  volume  = {87},
  pages   = {51--61},
  year    = {2005},
  doi     = {10.4064/ap87-0-5}
}

\end{document}